\documentclass[11pt, a4paper]{article}
\usepackage{amsmath}
\usepackage{wasysym}
\usepackage{latexsym}
\usepackage{amsfonts}
\usepackage{mathrsfs}
\usepackage{amssymb}
\usepackage{ifsym}
\usepackage{dsfont}
\usepackage[all]{xy}
\usepackage{authblk}

\newtheorem{Definitions1}{Definition}[section]

\newtheorem{Theorems1}{Theorem}[section]

\newtheorem{Coroll1}[Theorems1]{Corollary}
\newtheorem{Lemma1}[Theorems1]{Lemma}

\newtheorem{Examp1}{Example}[section]

\newenvironment{proof}[1][Proof]{\begin{trivlist}
\item[\hskip \labelsep {\bfseries #1}]}{\end{trivlist}}

\newcommand{\qed}{\nobreak \ifvmode \relax \else
      \ifdim\lastskip<1.5em \hskip-\lastskip
      \hskip1.5em plus0em minus0.5em \fi \nobreak
      \vrule height0.75em width0.5em depth0.25em\fi}

\begin{document}
\title{Decidable fragments of the Simple Theory of Types with Infinity and $\mathrm{NF}$ \thanks{The research of the first and the third author was supported in part by EPSRC grant EP/H026835.}}
\author[1]{Anuj Dawar}
\author[2]{Thomas Forster} 
\author[1]{Zachiri McKenzie}
\affil[1]{University of Cambridge Computer Laboratory\\ \texttt{firstname.lastname@cl.cam.ac.uk}}
\affil[2]{DPMMS, University of Cambridge\\
\texttt{T.Forster@dpmms.cam.ac.uk}}
\maketitle

\begin{abstract}
We identify complete fragments of the Simple Theory of Types with Infinity ($\mathrm{TSTI}$) and Quine's $\mathrm{NF}$ set theory. We show that $\mathrm{TSTI}$ decides every sentence $\phi$ in the language of type theory that is in one of the following forms:
\begin{itemize}
\item[(A)] $\phi= \forall x_1^{r_1} \cdots \forall x_k^{r_k} \exists y_1^{s_1} \cdots \exists y_l^{s_l} \theta$ where the superscripts denote the types of the variables, $s_1 > \ldots > s_l$ and $\theta$ is quantifier-free, 
\item[(B)] $\phi= \forall x_1^{r_1} \cdots \forall x_k^{r_k} \exists y_1^{s} \cdots \exists y_l^{s} \theta$ where the superscripts denote the types of the variables and $\theta$ is quantifier-free.
\end{itemize}
This shows that $\mathrm{NF}$ decides every stratified sentence $\phi$ in the language of set theory that is in one of the following forms:
\begin{itemize}
\item[(A')] $\phi= \forall x_1 \cdots \forall x_k \exists y_1 \cdots \exists y_l \theta$ where $\theta$ is quantifier-free and $\phi$ admits a stratification that assigns distinct values to all of the variable $y_1, \ldots, y_l$, 
\item[(B')] $\phi= \forall x_1 \cdots \forall x_k \exists y_1 \cdots \exists y_l \theta$ where $\theta$ is quantifier-free and $\phi$ admits a stratification that assigns the same value to all of the variables $y_1, \ldots, y_l$.
\end{itemize}  
\end{abstract}

\section[Introduction]{Introduction}

\indent Roland Hinnion showed in his thesis \cite{hin75} that {\sl Every consistent $\exists^*$ sentence in the language of set theory is a theorem of $\mathrm{NF}$} or, equivalently: {\sl Every finite binary structure can be embedded in every model of $\mathrm{NF}$}. Both these formulations invite generalisations.  On the one hand we find results like {\sl every countable binary structure can be embedded in every model of $\mathrm{NF}$} (this is theorem 4 of \cite{for87}) and on the other we can ask about the status of sentences with more quantifiers: $\forall^*\exists^*$ sentences in the first instance; it is the second that will be our concern here.\\
\\
\indent It is elementary to check that $\mathrm{NF}$ does not decide all $\forall^*\exists^*$ sentences, since the existence of Quine atoms ($x = \{x\}$) is consistent with, and independent of, $\mathrm{NF}$. However `$(\forall x)(x \not= \{x\})$' is not stratified, and this invites the conjecture that (i) $\mathrm{NF}$ decides all  stratified $\forall^*\exists^*$ sentences and that (ii) all  unstratified $\forall^*\exists^*$ sentences can be proved both  relatively consistent and independent by means of Rieger-Bernays permutation methods. It's with limb (i) of this conjecture that we are concerned here.\\
\\ 
\indent The foregoing is all about $\mathrm{NF}$; the connection with the Simple Theory of Types with Infinity ($\mathrm{TSTI}$) arises because of work of Ernst Specker \cite{spe62} and \cite{spe53}: $\mathrm{NF}$ decides all stratified $\forall^*\exists^*$ sentences of the language of set theory if and only if $\mathrm{TSTI} + \mathrm{Ambiguity}$ decides all $\forall^*\exists^*$ sentences of the language of type theory.
 
\begin{quote}{\bf Conjecture}: All models of $\mathrm{TSTI}$ agree on all $\forall^*\exists^*$ sentences.\end{quote}
 
\noindent It is towards a proof of this conjecture that our efforts in this paper are directed.\\
\\ 
\indent Observe that {\sl there is a total order of $V$} is consistent with and independent of $\mathrm{TST}$ and it can be said with three blocks of quantifiers: $$(\exists O)[(\forall x y \in O)(x \subseteq y \lor y \subseteq x) \land (\forall u v)( u \not= v \to (\exists x \in O)( u \in x \iff v \not\in x))]$$ making it $\exists^1\forall^6\exists^1$.

\section[Background and definitions]{Background and definitions} \label{Sec:Background}

The Simple Theory of Types is the simplification of the Ramified Theory of Types, the underlying system of \cite{rw08}, that was independently discovered by Frank Ramsey and Leon Chwistek. Following \cite{mat01} we use $\mathrm{TSTI}$ and $\mathrm{TST}$ to abbreviate the Simple Theory of Types with and without an axiom of infinity respectively. These theories are naturally axiomatised in a many-sorted language with sorts for each $n \in \mathbb{N}$.

\begin{Definitions1}
We use $\mathcal{L}_{\mathrm{TST}}$ to denote the $\mathbb{N}$-sorted language endowed with binary relation symbols $\in_n$ for each sort $n \in \mathbb{N}$. There are variables $x^n, y^n, z^n, \ldots$ for each sort $n \in \mathbb{N}$ and well-formed $\mathcal{L}_{\mathrm{TST}}$-formulae are built-up inductively from atomic formulae in the form $x^n \in_n y^{n+1}$ and $x^n = y^n$ using the connectives quantifiers of first-order logic.  
\end{Definitions1}

\noindent We refer to sorts of $\mathcal{L}_{\mathrm{TST}}$ as types. We will attempt to stick to the convention of denoting $\mathcal{L}_{\mathrm{TST}}$-structures using calligraphy letters ($\mathcal{M}, \mathcal{N}, \ldots$). A $\mathcal{L}_{\mathrm{TST}}$-structure $\mathcal{M}$ consists of domains $M_n$ for each type $n \in \mathbb{N}$ and interpretations of the relations $\in_n^{\mathcal{M}} \subseteq M_n \times M_{n+1}$ for each type $n \in \mathbb{N}$; we write $\mathcal{M}= \langle M_0, M_1, \ldots, \in_0^{\mathcal{M}}, \in_1^{\mathcal{M}}, \ldots \rangle$. If \mbox{$\mathcal{M}= \langle M_0, M_1, \ldots, \in_0^{\mathcal{M}}, \in_1^{\mathcal{M}}, \ldots \rangle$} is an $\mathcal{L}_{\mathrm{TST}}$-structure then we call the elements of $M_0$ atoms.

\begin{Definitions1}
We use $\mathrm{TST}$ to denote the $\mathcal{L}_{\mathrm{TST}}$-theory with axioms
\begin{itemize}
\item[](Extensionality) for all $n \in \mathbb{N}$, 
$$\forall x^{n+1} \forall y^{n+1} (x^{n+1}= y^{n+1} \iff \forall z^n(z^{n} \in_n x^{n+1} \iff z^n \in y^{n+1})),$$  
\item[](Comprehension) for all $n \in \mathcal{N}$ and for all well-formed $\mathcal{L}_{\mathrm{TST}}$-formulae $\phi(x^n, \vec{z})$,
$$\forall \vec{z} \exists y^{n+1} \forall x^n (x^n \in_n y^{n+1} \iff \phi(x^n, \vec{z})).$$  
\end{itemize} 
\end{Definitions1}

Comprehension ensures that every successor type is closed under the set-theoretic operations: union ($\cup$), intersection ($\cap$), difference ($\backslash$) and symmetric difference ($\triangle$). For all $n \in \mathbb{N}$, we use $\emptyset^{n+1}$ to denote the point at type $n+1$ which contains no points from type $n$ and we use $V^{n+1}$ to denote the point at type $n+1$ that contains every point from type $n$. The Wiener-Kuratowski ordered pair allows us to code ordered pairs in the form $\langle x, y \rangle$ as objects in $\mathrm{TST}$ which have type two higher than the type of $x$ and $y$. Functions, as usual, are thought of as collections of ordered pairs. This means that a function $f: X \longrightarrow Y$ will be coded by an object in $\mathrm{TST}$ that has type two higher than the type of $X$ and $Y$. The theory $\mathrm{TSTI}$ is obtained from $\mathrm{TST}$ by asserting the existence of a Dedekind infinite collection at type $1$.

\begin{Definitions1}
We use $\mathrm{TSTI}$ to denote the $\mathcal{L}_{\mathrm{TST}}$-theory obtained from $\mathrm{TST}$ by adding the axiom
$$\exists x^1 \exists f^3(f^3: x^1 \longrightarrow x^1 \textrm{ is injective but not surjective}).$$ 
\end{Definitions1}

Let $X$ be a set. If the $\mathcal{L}_{\mathrm{TST}}$-structure $\mathcal{M}= \langle M_0, M_1, \ldots, \in_0, \in_1, \ldots \rangle$ is defined by $M_n= \mathcal{P}^n(X)$ and $\in_n^\mathcal{M}= \in \upharpoonright \mathcal{P}^n(X) \times \mathcal{P}^{n+1}(X)$ for all $n \in \mathbb{N}$, then $\mathcal{M} \models \mathrm{TST}$. If $m \in \mathbb{N}$ and $|X|= m$ then $\mathcal{M}$ is the unique, up to isomorphism, model of $\mathrm{TST}$ with exactly $m$ atoms and we say that $\mathcal{M}$ is finitely generated by $m$ atoms. Alternatively, if $X$ is Dedekind infinite then $\mathcal{M} \models \mathrm{TSTI}$. This shows that $\mathrm{ZFC}$ proves the consistency of $\mathrm{TSTI}$. In fact, in \cite{mat01} it is shown that $\mathrm{TSTI}$ is equiconsistent with Mac Lane Set Theory.\\
\\
\indent We say that an $\mathcal{L}^\prime$-theory $T$ decides an $\mathcal{L}^\prime$-sentence $\phi$ if and only if $T \vdash \phi$ or $T \vdash \neg \phi$. The Completeness Theorem implies that $T$ decides $\phi$ if and only if $\phi$ holds in all $\mathcal{L}^\prime$-structures $\mathcal{M} \models T$, or $\neg \phi$ holds in all $\mathcal{L}^\prime$-structures $\mathcal{M} \models T$.    

\begin{Definitions1}
We say that a $\mathcal{L}_{\mathrm{TST}}$-sentence $\phi$ is $\exists^* \forall^*$ if and only if\\ 
\mbox{$\phi= \exists x_1^{r_1} \cdots \exists x_k^{r_k} \forall y_1^{s_1} \cdots \forall y_l^{s_l} \theta$} where $\theta$ is quantifier-free.
\end{Definitions1}        

\begin{Definitions1}
We say that an $\mathcal{L}_{\mathrm{TST}}$-sentence $\phi$ is $\forall^* \exists^*$ if and only if\\ 
\mbox{$\phi= \forall x_1^{r_1} \cdots \forall x_k^{r_k} \exists y_1^{s_1} \cdots \exists y_l^{s_l} \theta$} where $\theta$ is quantifier-free.
\end{Definitions1}

We will show that $\mathrm{TSTI}$ decides a significant fragment of the $\forall^* \exists^*$ sentences (and thus it also decides the $\exists^* \forall^*$ sentences that are logically equivalent to the negation of these $\forall^* \exists^*$ sentences). We achieve this result by showing that every sentence or negation of a sentence in this fragment that is true in some model of $\mathrm{TSTI}$ is true in all models of $\mathrm{TST}$ that are finitely generated by sufficiently many atoms.

\begin{Definitions1}
We say that an $\mathcal{L}_{\mathrm{TST}}$-sentence $\phi$ has the finitely generated model property if and only if, if there exists an $\mathcal{N} \models \mathrm{TSTI}+\phi$ then there exists a $k \in \mathbb{N}$ such that for all $m \geq k$, if $\mathcal{M} \models \mathrm{TST}$ is finitely generated by $m$ atoms then $\mathcal{M} \models \phi$. 
\end{Definitions1}

\noindent Note that if $\Gamma$ is class of $\mathcal{L}_{\mathrm{TST}}$-sentences that have the finitely generated model property and $\Gamma$ is closed under negations then $\mathrm{TST}$ decides every sentence in $\Gamma$.\\
\\
\indent In \cite{qui37} Willard van Orman Quine introduces a set theory by identifying a syntactic condition on formulae in the single sorted language of set theory that captures the restricted comprehension available in $\mathrm{TST}$. This set theory has been dubbed `New Foundations' ($\mathrm{NF}$) after the title of \cite{qui37}. We will use $\mathcal{L}$ to denote the language of set theory --- the language of first-order logic endowed with a binary relation symbol $\in$ whose intended interpretation is membership. Before giving the axioms of $\mathrm{NF}$ we first recall Quine's definition of a stratified formulae. If $\phi$ is an $\mathcal{L}$-formula then we use $\mathbf{Var}(\phi)$ to denote the set of variables (both free and bound) which appear in $\phi$. 

\begin{Definitions1}
Let $\phi(x_1, \ldots, x_n)$ be an $\mathcal{L}$-formula. We say that $\sigma: \mathbf{Var}(\phi) \longrightarrow \mathbb{N}$ is a stratification of $\phi$ if and only if  
\begin{itemize}
\item[(i)] if `$x \in y$' is a subformula of $\phi$ then $\sigma(\textrm{`}y\textrm{'})= \sigma(\textrm{`}x\textrm{'})+1$, 
\item[(ii)] if `$x = y$' is a subformula of $\phi$ then $\sigma(\textrm{`}y\textrm{'})= \sigma(\textrm{`}x\textrm{'})$.
\end{itemize}
If there exists a stratification of $\phi$ then we say that $\phi$ is stratified.
\end{Definitions1}

Let $\phi$ be an $\mathcal{L}$-formula. Note that $\sigma: \mathbf{Var}(\phi) \longrightarrow \mathbb{N}$ is a stratification of $\phi$ if and only if the formula obtained by decorating every variable appearing in $\phi$ with the type given by $\sigma$ yields a well-formed $\mathcal{L}_{\mathrm{TST}}$-formula. Conversely, let $\theta$ be a well-formed $\mathcal{L}_{\mathrm{TST}}$-formula and let $\phi$ an $\mathcal{L}$-formula obtained for $\theta$ by deleting the types from the variables appearing in $\theta$ while ensuring (by relabeling variables) that no two distinct variables in $\theta$ become the same variable in $\phi$. Then the $\mathcal{L}$-formula $\phi$ is stratified and the function which sends a variable in $\phi$ to the type index of the corresponding variable in $\theta$ is a stratification.   

\begin{Definitions1}
Let $\phi$ be an $\mathcal{L}$-formula with stratification $\sigma: \mathbf{Var}(\phi)$. We use $\phi^{(\sigma)}$ to denote the $\mathcal{L}_{\mathrm{TST}}$-formula obtained by assigning each variable `$x$' appearing $\phi$ the type $\sigma(\textrm{`}x\textrm{'})$. 
\end{Definitions1}

$\mathrm{NF}$ is the $\mathcal{L}$-theory with the axiom of extensionality and comprehension for all stratified $\mathcal{L}$-formulae.

\begin{Definitions1}
We use $\mathrm{NF}$ to denote the $\mathcal{L}$-theory with axioms
\begin{itemize}
\item[](Extensionality) $\forall x \forall y (x=y \iff \forall z(z \in x \iff z \in y))$,  
\item[](Stratified Comprehension) for all stratified $\phi(x, \vec{z})$,
$$\forall \vec{z} \exists y \forall x (x \in y \iff \phi(x, \vec{z})).$$ 
\end{itemize}
\end{Definitions1}

We direct the interested reader to \cite{for95} for detailed treatment of $\mathrm{NF}$. One interesting feature of $\mathrm{NF}$ is that it refutes the Axiom of Choice and so proves the Axiom of Infinity (see \cite{spe53}). There is a strong connection between the theories $\mathrm{NF}$ and $\mathrm{TSTI}$. \cite{spe62} shows that models of $\mathrm{NF}$ can be obtained from models of $\mathrm{TSTI}$ plus the scheme $\phi \iff \phi^+$, for all $\mathcal{L}_{\mathrm{TST}}$-sentences $\phi$, where $\phi^+$ is obtained from $\phi$ by incrementing the types of all the variables appearing in $\phi$. Conversely, let $\mathcal{M}= \langle M, \in^{\mathcal{M}} \rangle$ be an $\mathcal{L}$-structure with $\mathcal{M} \models \mathrm{NF}$. The $\mathcal{L}_{\mathrm{TST}}$-structure $\mathcal{N}= \langle N_0, N_1, \ldots, \in_0^{\mathcal{N}}, \in_1^{\mathcal{N}}, \ldots \rangle$ defined by $N_n= M$ and $\in_n^{\mathcal{N}}= \in^{\mathcal{M}}$ is such that $\mathcal{N} \models \mathrm{TSTI}$. Moreover, if $\phi$ is an $\mathcal{L}$-sentence with stratification $\sigma: \mathbf{Var}(\phi) \longrightarrow \mathbb{N}$ and $\mathcal{M} \models \phi$ then $\mathcal{N} \models \phi^{(\sigma)}$. This immediately shows that a decidable fragment of $\mathrm{TSTI}$ yields a decidable fragment of $\mathrm{NF}$.

\begin{Theorems1} \label{Th:DecidableFragmentsOfNF}
Let $\phi$ be an $\mathcal{L}$-sentence with stratification $\sigma: \mathbf{Var}(\phi) \longrightarrow \mathbb{N}$. If $\mathrm{TSTI}$ decides $\phi^{(\sigma)}$ then $\mathrm{NF}$ decides $\phi$.
\Square
\end{Theorems1}

\section[$\exists^* \forall^*$ sentences have the finitely generated model property]{$\exists^* \forall^*$ sentences have the finitely generated model property}

In this section we prove that all $\exists^* \forall^*$ sentences have the finitely generated model property. This result follows from the fact that if $\mathcal{N}$ is a model of $\mathrm{TSTI}$, $a_1^{r_1}, \ldots, a_k^{r_k} \in \mathcal{N}$ with $r_1 \leq \ldots \leq r_k$ and $\mathcal{M}$ is a model of $\mathrm{TST}$ that is finitely generated by sufficiently many atoms then there is an embedding of $\mathcal{M}$ into $\mathcal{N}$ with $a_1^{r_1}, \ldots, a_k^{r_k}$ in the range. Given $k \in \mathbb{N}$ we define the function $\mathbf{G}_k: \mathbb{N} \longrightarrow \mathbb{N}$ by recursion
\begin{equation}
\mathbf{G}_k(0)= k \textrm{ and } \mathbf{G}_k(n+1)= \binom{\mathbf{G}_k(n)}{2}+ k. 
\end{equation}

\begin{Lemma1} \label{Th:EmeddingProperty}
Let $\mathcal{N} \models \mathrm{TSTI}$ and let \mbox{$a_1^{r_1}, \ldots, a_k^{r_k} \in \mathcal{N}$} with $r_1 \leq \ldots \leq r_k$. If $\mathcal{M} \models \mathrm{TST}$ is finitely generated by at least $\mathbf{G}_k(r_k)$ atoms then there exists a sequence $\langle f_n \mid n \in \mathbb{N} \rangle$ such that for all $n \in \mathbb{N}$,
\begin{itemize}
\item[(i)] $f_n: M_n \longrightarrow N_n$ is injective, 
\item[(ii)] for all $x \in M_n$ and for all $y \in M_{n+1}$,
$$\mathcal{M} \models x \in_n y \textrm{ if and only if } \mathcal{N} \models f_n(x) \in_n f_{n+1}(y),$$ 
\item[(iii)] $$a_1^{r_1}, \ldots, a_k^{r_k} \in \bigcup_{m \in \mathbb{N}} \mathrm{rng}(f_m).$$ 
\end{itemize} 
\end{Lemma1}

\begin{proof}
Let $\mathcal{N}= \langle N_0, N_1, \ldots, \in_0^{\mathcal{N}}, \in_1^{\mathcal{N}}, \ldots \rangle$ be such that $\mathcal{N} \models \mathrm{TSTI}$ and let $a_1^{r_1}, \ldots, a_k^{r_k} \in \mathcal{N}$ with $r_1 \leq \ldots \leq r_k$. Let $\mathcal{M}= \langle M_0, M_1, \ldots, \in_0^{\mathcal{M}}, \in_1^{\mathcal{M}}, \ldots \rangle$ be such that $\mathcal{M} \models \mathrm{TST}$ is finitely generated and $|M_0| \geq \mathbf{G}_k(r_k)$. We begin by defining $C \subseteq \mathcal{N}$ such that $|C \cap N_0| \leq \mathbf{G}_k(r_k)$ and for any two points $x \neq y$ in $C$ that are not atoms, there exists a point $z$ in $C$ which $\mathcal{N}$ believes is in the symmetric difference of $x$ and $y$. Define $C_0= \{a_1^{r_1}, \ldots, a_k^{r_k} \} \subseteq \mathcal{N}$. Note that $|C_0 \cap N_{r_k}| \leq \mathbf{G}_k(0)= k$ and for all $0\leq m < r_k$, $|C_0 \cap N_m| \leq k$. For $0 < n \leq r_k$ we recursively define $C_n \subseteq \mathcal{N}$ which satisfies
\begin{itemize}
\item[(I)] $|C_n \cap N_{r_k-n}| \leq \mathbf{G}_k(n)$,
\item[(II)] for all $0 \leq m < r_k-n$, $|C_n \cap N_m| \leq k$. 
\end{itemize}
Suppose that $n < r_k$ and $C_n \subseteq \mathcal{N}$ has been defined and satisfies (I) and (II). For all $y, z \in N_{r_k-n}$ with $y \neq z$, let $\gamma_{\{y, z\}} \in N_{r_k-(n+1)}$ be such that
$$\mathcal{N} \models \gamma_{\{y, z\}} \in_{r_k-(n+1)} y \triangle z.$$
Define
$$C_{n+1}= C_n \cup \{ \gamma_{\{y, z\}} \mid \{y, z\} \in [N_{r_k-n} \cap C_n]^2 \}.$$
It follows from (I) and (II) that
$$|C_{n+1} \cap N_{r_k-(n+1)}| \leq |C_n \cap N_{r_k-(n+1)}|+\binom{|C_n \cap N_{r_k-n}|}{2} \leq k + \binom{\mathbf{G}_k(n)}{2}= \mathbf{G}_k(n+1)$$
and for all $0 \leq m < r_k-(n+1)$, $|C_{n+1} \cap N_m| \leq k$. Now, let $C= C_{r_k}$. This recursion ensures that $|C \cap N_0| \leq \mathbf{G}_k(r_k)$.\\
We now turn to defining the family of maps $\langle f_n \mid n \in \mathbb{N} \rangle$ which embed $\mathcal{M}$ into $\mathcal{N}$. We define the sequence $\langle f_n \mid n \in \mathbb{N} \rangle$ by induction. Let $C^\prime= C \cap N_0$. Let $f_0: M_0 \longrightarrow N_0$ be an injection such that $C^\prime \subseteq \mathrm{rng}(f_0)$. Suppose that $\langle f_0, \ldots, f_n \rangle$ has been defined such that
\begin{itemize}
\item[(I')] for all $0 \leq j \leq n$, $f_j: M_j \longrightarrow N_j$ is injective,
\item[(II')] for all $0 \leq j < n$, for all $x \in M_j$ and for all $y \in M_{j+1}$,
$$\mathcal{M} \models x \in_j y \textrm{ if and only if } \mathcal{N} \models f_j(x) \in_j f_{j+1}(y),$$
\item[(III')] for all $0 \leq j \leq n$, $C \cap N_j \subseteq \mathrm{rng}(f_j)$.
\end{itemize}
If $0 \leq j \leq n$ and $x \in M_{j+1}$ then we use $f_j``x$ to denote the point in $N_{j+1}$ such that $\mathcal{N} \models f_j``x= \{ f_j(y) \mid \mathcal{M} \models y \in_j x \}$. Note that, since $\mathcal{M}$ is finitely generated, for all $x \in M_{j+1}$, $f_j``x$ exists in $\mathcal{N}$. We define $f_{n+1}: M_{n+1} \longrightarrow N_{n+1}$ by
$$f_{n+1}(x)= \left\{ \begin{array}{ll}
\gamma & \textrm{if } \gamma \in C \cap N_{n+1} \textrm{ and } \mathcal{N} \models f_n``x= \gamma \cap f_n``(V^{n+1})^\mathcal{M}\\
f_n``x & \textrm{otherwise}
\end{array} \right.$$
We first need to show that the map $f_{n+1}$ is well-defined. Suppose that $\xi_1, \xi_2 \in C \cap N_{n+1}$ with $\xi_1 \neq \xi_2$ and $x \in M_{n+1}$ are such that
$$\mathcal{N} \models f_n``x= \xi_1 \cap f_n``(V^{n+1})^\mathcal{M} \textrm{ and } \mathcal{N} \models f_n``x= \xi_2 \cap f_n``(V^{n+1})^\mathcal{M}.$$
Now, there is a $\gamma \in C \cap N_n$ such that $\mathcal{N} \models \gamma \in_n \xi_1 \triangle \xi_2$. By (III'), $\gamma \in \mathrm{rng}(f_n)$, which is a contradiction. Therefore $f_{n+1}$ is well-defined.\\
The fact that $f_n$ is injective ensures that $f_{n+1}$ is injective.\\
We now turn to showing that the sequence $\langle f_0, \ldots, f_{n+1} \rangle$ satisfies (II'). Let $x \in M_n$ and let $y \in M_{n+1}$. There are two cases. Firstly, suppose that $f_{n+1}(y)= \gamma \in C$. Therefore $\mathcal{N} \models f_n``y= \gamma \cap f``(V^{n+1})^\mathcal{M}$. If $\mathcal{M} \models x \in_n y$ then $\mathcal{N}\models f_n(x) \in_n f_n``y$ and so $\mathcal{N} \models f_n(x) \in_n f_{n+1}(y)$. Conversely, if $\mathcal{N} \models f_n(x) \in_n \gamma$ then $\mathcal{N} \models f_n(x) \in_n f_n``y$ and so $\mathcal{M} \models x \in_n y$. The second case is when $f_{n+1}(y)= f_n``y$. In this case it is clear that
$$\mathcal{M} \models x \in_n y \textrm{ if and only if } \mathcal{N} \models f_n(x) \in_n f_{n+1}(y).$$
This shows that the sequence $\langle f_0, \ldots, f_{n+1} \rangle$ satisfies (II').\\
This concludes the induction step of the construction and shows that we can construct a sequence $\langle f_n \mid n \in \mathbb{N} \rangle$ that satisfies (i)-(iii).     
\Square
\end{proof}

This embedding property allows us to show that every $\exists^* \forall^*$ sentence has the finitely generated model property.

\begin{Theorems1} \label{Th:ExistentialUniversalSentencesHaveFinitelyGeneratedModelProperty}
Let $\phi= \exists x_1^{r_1} \cdots \exists x_k^{r_k} \forall y_1^{s_1} \cdots \forall y_l^{s_l} \theta$ where $r_1 \leq \ldots \leq r_k$ and $\theta$ is quantifier-free. If $\mathcal{N} \models \mathrm{TSTI}+\phi$ and $\mathcal{M} \models \mathrm{TST}$ is finitely generated by at least $\mathbf{G}_k(r_k)$ atoms then $\mathcal{M} \models \phi$.
\end{Theorems1}

\begin{proof}
Let $\mathcal{N}= \langle N_0, N_1, \ldots, \in_0^{\mathcal{N}}, \in_1^{\mathcal{N}}, \ldots \rangle$ be such that $\mathcal{N} \models \mathrm{TSTI}+\phi$. Let $\mathcal{M}= \langle M_0, M_1, \ldots, \in_0^{\mathcal{M}}, \in_1^{\mathcal{M}}, \ldots \rangle$ be such that $\mathcal{M} \models \mathrm{TST}$ and $\mathcal{M}$ is finitely generated by at least $\mathbf{G}_k(r_k)$ atoms. Let $a_1^{r_1}, \ldots, a_k^{r_k} \in \mathcal{N}$ be such that
$$\mathcal{N} \models \forall y_1^{s_1} \cdots \forall y_l^{s_l} \theta[a_1^{r_1}, \ldots, a_k^{r_k}].$$
Using Lemma \ref{Th:EmeddingProperty} we can find a sequence $\langle f_n \mid n \in \mathbb{N} \rangle$ such that
\begin{itemize}
\item[(i)] $f_n: M_n \longrightarrow N_n$ is injective, 
\item[(ii)] for all $x \in M_n$ and for all $y \in M_{n+1}$,
$$\mathcal{M} \models x \in_n y \textrm{ if and only if } \mathcal{N} \models f_n(x) \in f_{n+1}(y),$$ 
\item[(iii)] $$a_1^{r_1}, \ldots, a_k^{r_k} \in \bigcup_{m \in \mathbb{N}} \mathrm{rng}(f_m).$$ 
\end{itemize}
Let $b_1^{r_1}, \ldots, b_k^{r_k} \in \mathcal{M}$ be such that for all $1 \leq j \leq k$, $f_{r_j}(b_j^{r_j})= a_j^{r_j}$. Let $c_1^{s_1}, \ldots, c_l^{s_l} \in \mathcal{M}$. Since $\mathcal{N} \models \theta[a_1^{r_1}, \ldots, a_k^{r_k}, f_{s_1}(c_1^{s_1}), \ldots, f_{s_l}(c_l^{s_l})]$, it follows that
$$\mathcal{M} \models \theta[b_1^{r_1}, \ldots, b_k^{r_k}, c_1^{s_1}, \ldots, c_l^{s_l}].$$
Therefore
$$\mathcal{M} \models \forall y_1^{s_1} \cdots \forall y_l^{s_l}\theta[b_1^{r_1}, \ldots, b_k^{r_k}],$$
which proves the theorem.
\Square
\end{proof}

\section[Decidable fragments of the $\forall^* \exists^*$ sentences]{Decidable fragments of the $\forall^* \exists^*$ sentences}

In this section we will show that $\mathrm{TSTI}$ decides every $\forall^* \exists^*$ sentence $\phi$ that is in one of the following forms:
\begin{itemize}
\item[(A)] $\phi= \forall x_1^{r_1} \cdots \forall x_k^{r_k} \exists y_1^{s_1} \cdots \exists y_l^{s_l} \theta$ where $s_1 > \ldots > s_l$ and $\theta$ is quantifier-free, 
\item[(B)] $\phi= \forall x_1^{r_1} \cdots \forall x_k^{r_k} \exists y_1^{s} \cdots \exists y_l^{s} \theta$ where $\theta$ is quantifier-free.
\end{itemize}
By applying Theorem \ref{Th:DecidableFragmentsOfNF} it then follows that $\mathrm{NF}$ decides every stratified $\mathcal{L}$-sentence $\phi$ that is in one of the following forms:
\begin{itemize}
\item[(A')] $\phi= \forall x_1 \cdots \forall x_k \exists y_1 \cdots \exists y_l \theta$ where $\theta$ is quantifier-free and $\sigma: \mathbf{Var}(\phi) \longrightarrow \mathbb{N}$ is a stratification of $\phi$ that assigns distinct values to all of the variables $y_1, \ldots, y_l$, 
\item[(B')] $\phi= \forall x_1 \cdots \forall x_k \exists y_1 \cdots \exists y_l \theta$ where $\theta$ is quantifier-free and $\sigma: \mathbf{Var}(\phi) \longrightarrow \mathbb{N}$ is a stratification of $\phi$ that assigns the same value to all of the variables $y_1, \ldots, y_l$.
\end{itemize}
Throughout this section we will fix $k, l \in \mathbb{N}$ and a sequence $r_1 \leq \ldots \leq r_k$ that will represent the types of the universally quantified variables in a $\forall^* \exists^*$ sentence. Let $k^\prime$ be the number of distinct elements in the list $r_1, \ldots, r_k$. Let $K_1, \ldots, K_{k^\prime}$ be the multiplicities of the elements in the list $r_1, \ldots, r_k$, so $k= \sum_{1 \leq i \leq k^\prime} K_i$, and let $K= \max\{ K_1, \ldots, K_{k^\prime}, l\}$. We also fix structures $\mathcal{N}= \langle N_0, N_1, \ldots, \in_0^{\mathcal{N}}, \in_1^{\mathcal{N}}, \ldots \rangle$ with $\mathcal{N} \models \mathrm{TSTI}$ and $\mathcal{M}= \langle M_0, M_1, \ldots, \in_0^{\mathcal{M}}, \in_1^{\mathcal{M}}, \ldots \rangle$ with $\mathcal{M} \models \mathrm{TST}$ finitely generated by at least $(2^K)^{k^\prime+2}$ atoms. Let $a_1^{r_1}, \ldots, a_k^{r_k} \in \mathcal{M}$.\\
\\
\indent Our approach will be to define colour classes $\mathcal{C}_{i, j}$, the elements of which we will call colours, and functions $c_{i, j}^{\mathcal{M}}: M_i \longrightarrow \mathcal{C}_{i, j}$ and $c_{i, j}^{\mathcal{N}}: N_i \longrightarrow \mathcal{C}_{i, j}$, which we will call colourings, for all $i \in \mathbb{N}$ and for all $0 \leq j \leq k^\prime$. For all $0 < j \leq k^\prime$, the colourings $c_{i, j}^{\mathcal{M}}$ will be defined using the elements $a_1^{r_1}, \ldots, a_{j^\prime}^{r_{j^\prime}}$ where $j^\prime= \sum_{1 \leq m \leq j} K_m$, and in the process of defining the colourings $c_{i, j}^{\mathcal{N}}$ we will construct corresponding elements $b_1^{r_1}, \ldots, b_{j^\prime}^{r_{j^\prime}} \in \mathcal{N}$. The colourings will be designed with the following properties:
\begin{itemize}
\item[(i)] For a fixed colour $\alpha$ in some $\mathcal{C}_{i, j}$, the property of being an element of $\mathcal{N}$ that is given colour $\alpha$ by $c_{i, j}^{\mathcal{N}}$ will be definable by an $\mathcal{L}_{\mathrm{TST}}$-formula, $\Phi_{i, j, \alpha}$, with parameters over $\mathcal{N}$.  
\item[(ii)] The colour given to an element $x$ in $\mathcal{M}$ (or $\mathcal{N}$) by the colouring $c_{i, j}^{\mathcal{M}}$ (respectively $c_{i, j}^{\mathcal{N}}$) will tell us which quantifier-free $\mathcal{L}_{\mathrm{TST}}$-formulae with parameters $a_1^{r_1}, \ldots, a_{j^\prime}^{r_{j^\prime}}$ (respectively $b_1^{r_1}, \ldots, b_{j^\prime}^{r_{j^\prime}}$), where $j^\prime= \sum_{1 \leq m \leq j} K_m$, are satisfied by $x$ in $\mathcal{M}$ (respectively $\mathcal{N}$).   
\item[(iii)] For every colour $\beta$ in $\mathcal{C}_{i, j}$, the colour given to an element $x$ in $\mathcal{M}$ (or $\mathcal{N}$) by the colouring $c_{i+1, j}^{\mathcal{M}}$ (respectively $c_{i+1, j}^{\mathcal{N}}$) will tell us whether or not there is an element $y$ in $\mathcal{M}$ (respectively $\mathcal{N}$) such that $\mathcal{M} \models y \in_i x$ (respectively $\mathcal{N} \models y \in_i x$) and $y$ is given colour $\beta$ by $c_{i, j}^{\mathcal{M}}$ (respectively $c_{i, j}^{\mathcal{N}}$).
\item[(iv)] For every colour $\beta$ in $\mathcal{C}_{i, j}$, the colour given to an element $x$ in $\mathcal{M}$ (or $\mathcal{N}$) by the colouring $c_{i+1, j}^{\mathcal{M}}$ (respectively $c_{i+1, j}^{\mathcal{N}}$) will tell us whether or not there is an element $y$ in $\mathcal{M}$ (respectively $\mathcal{N}$) such that $\mathcal{M} \models y \notin_i x$ (respectively $\mathcal{N} \models y \notin_i x$) and $y$ is given colour $\beta$ by $c_{i, j}^{\mathcal{M}}$ (respectively $c_{i, j}^{\mathcal{N}}$). 
\end{itemize}
Note that since $\mathcal{M}$ is finitely generated, the analogue of condition (i) automatically holds for $\mathcal{M}$.\\
\\
\indent Before defining the colour classes $\mathcal{C}_{i, j}$ and the colourings $c_{i,j}^{\mathcal{M}}$ and $c_{i, j}^{\mathcal{N}}$ we first introduce the following definitions:

\begin{Definitions1}
Let $m \in \mathbb{N}$. We say that a colour $\alpha \in \mathcal{C}_{i, j}$ is $m$-special with respect to a colouring $f: X \longrightarrow \mathcal{C}_{i, j}$ if and only if 
$$|\{ x \in X \mid f(x)= \alpha\}|= m.$$
If $\alpha \in \mathcal{C}_{i, j}$ is $0$-special then we say that $\alpha$ is forbidden.
\end{Definitions1}

\begin{Definitions1}
Let $m \in \mathbb{N}$. We say that a colour $\alpha \in \mathcal{C}_{i, j}$ is $m$-abundant with respect to a colouring $f: X \longrightarrow \mathcal{C}_{i, j}$ if and only if 
$$|\{ x \in X \mid f(x)= \alpha\}|\geq m.$$
\end{Definitions1}

\begin{Definitions1}
Let $J \in \mathbb{N}$. We say that colourings $f: X \longrightarrow \mathcal{C}_{i, j}$ and $g: Y \longrightarrow \mathcal{C}_{i, j}$ are $J$-similar if and only if for all $0 \leq m < J$ and for all $\alpha \in \mathcal{C}_{i, j}$,
$$\alpha \textrm{ is } m \textrm{-special w.r.t. } f \textrm{ if and only if } \alpha \textrm{ is } m \textrm{-special w.r.t. } g.$$ 
\end{Definitions1}

The colour classes $\mathcal{C}_{i, j}$ and colourings $c_{i, j}^\mathcal{M}$ and $c_{i, j}^\mathcal{N}$ for all $i \in \mathbb{N}$ and for all $0 \leq j \leq k^\prime$ will be defined by a two-dimensional recursion. At each stage of the construction we will ensure that $c_{i, j}^\mathcal{M}$ and $c_{i, j}^\mathcal{N}$ are $(2^K)^{k^\prime-j+2}$-similar.\\
\\
Let $\mathcal{C}_{0, 0}= \{0\}$. Define $c_{0, 0}^{\mathcal{M}}: M_0 \longrightarrow \mathcal{C}_{0, 0}$ by
$$c_{0, 0}^{\mathcal{M}}(x)= 0 \textrm{ for all } x \in M_0.$$
Define $c_{0, 0}^{\mathcal{N}}: N_0 \longrightarrow \mathcal{C}_{0, 0}$ by
$$c_{0, 0}^{\mathcal{N}}(x)= 0 \textrm{ for all } x \in N_0.$$
Let $\Phi_{0, 0, 0}(x^0)$ be the $\mathcal{L}_{\mathrm{TST}}$-formula $x^0=x^0$. Note that for all $x \in N_0$,
$$\mathcal{N} \models \Phi_{0, 0, 0}[x] \textrm{ if and only if } c_{0, 0}^\mathcal{N}(x)= 0.$$

\begin{Lemma1}
The colourings $c_{0, 0}^{\mathcal{M}}$ and $c_{0, 0}^{\mathcal{N}}$ are $(2^K)^{k^\prime+2}$-similar. 
\end{Lemma1}

\begin{proof}
This follows immediately from the fact that $|M_0| \geq (2^K)^{k^\prime+2}$.
\Square
\end{proof}

We now turn to defining the colour classes $\mathcal{C}_{i, 0}$ and colourings $c_{i, 0}^{\mathcal{M}}: M_i \longrightarrow \mathcal{C}_{i, 0}$ and $c_{i, 0}^{\mathcal{N}}: N_i \longrightarrow \mathcal{C}_{i, 0}$ for all $i \in \mathbb{N}$. Suppose that we have defined the colour class $\mathcal{C}_{n, 0}$ with a canonical ordering, colourings $c_{n, 0}^{\mathcal{M}}: M_n \longrightarrow \mathcal{C}_{n, 0}$ and $c_{n, 0}^{\mathcal{N}}: N_i \longrightarrow \mathcal{C}_{n, 0}$ and $\mathcal{L}_{\mathrm{TST}}$-formulae $\Phi_{n, 0, \alpha}(x^n)$ for all $\alpha \in \mathcal{C}_{n, 0}$ with the following properties:
\begin{itemize}
\item[(I)] $c_{n, 0}^{\mathcal{M}}$ and $c_{n, 0}^{\mathcal{N}}$ are $(2^K)^{k^\prime+2}$-similar,
\item[(II)] for all $\alpha \in \mathcal{C}_{n, 0}$ and for all $x \in N_n$,
$$\mathcal{N} \models \Phi_{n, 0, \alpha}[x] \textrm{ if and only if } c_{n, 0}^{\mathcal{N}}(x)= \alpha.$$
\end{itemize}     
Let $\mathcal{C}_{n, 0}= \{ \alpha_1, \ldots, \alpha_q \}$ be the enumeration obtained from the canonical ordering. Define $\mathcal{C}_{n+1, 0}= 2^{2 \cdot q}$ --- the set of all 0-1 sequences of length $2 \cdot q$. Define $c_{n+1, 0}^{\mathcal{M}}: M_{n+1} \longrightarrow \mathcal{C}_{n+1, 0}$ such that for all $x \in M_{n+1}$,
$$c_{n+1, 0}^{\mathcal{M}}(x)= \langle f_1, \ldots, f_q, g_1, \ldots, g_q \rangle$$
$$\textrm{where } f_i= \left\{ \begin{array}{ll}
0 & \textrm{if for all } y \in M_n, \textrm{ if } c_{n, 0}^{\mathcal{M}}(y)= \alpha_i \textrm{ then } \mathcal{M}\models y \notin_n x\\
1 & \textrm{if there exists } y \in M_n, \textrm{ s.t. } c_{n, 0}^{\mathcal{M}}(y)= \alpha_i \textrm{ and } \mathcal{M} \models y \in_n x
\end{array} \right.$$
$$\textrm{and } g_i= \left\{ \begin{array}{ll}
0 & \textrm{if for all } y \in M_n, \textrm{ if } c_{n, 0}^{\mathcal{M}}(y)= \alpha_i \textrm{ then } \mathcal{M} \models y \in_n x\\
1 & \textrm{if there exists } y \in M_n \textrm{ s.t. } c_{n, 0}^{\mathcal{M}}(y)= \alpha_i \textrm{ and } \mathcal{M} \models y \notin_n x
\end{array} \right.$$

\begin{Examp1}
Using this definition we get $\mathcal{C}_{1, 0}= \{ \langle 0, 0 \rangle, \langle 1, 0 \rangle, \langle 0, 1 \rangle, \langle 1, 1 \rangle \}$. There are no $x \in M_1$ which are given the colour $\langle 0, 0 \rangle$ by $c_{1, 0}^{\mathcal{M}}$. The only point in $M_1$ which is given the colour $\langle 1, 0 \rangle$ by $c_{1, 0}^{\mathcal{M}}$ is $(V^1)^{\mathcal{M}}$. Similarly, the only point in $M_1$ which is given the colour $\langle 0, 1 \rangle$ by $c_{1, 0}^{\mathcal{M}}$ is $(\emptyset^1)^{\mathcal{M}}$. Every other point in $M_1$ is given the colour $\langle 1, 1 \rangle$ by $c_{1, 0}^{\mathcal{M}}$. 
\end{Examp1}

We define the colouring $c_{n+1, 0}^{\mathcal{N}}: N_{n+1} \longrightarrow \mathcal{C}_{n+1, 0}$ identically. Define $c_{n+1, 0}^{\mathcal{N}}: N_{n+1} \longrightarrow \mathcal{C}_{n+1, 0}$ such that for all $x \in N_{n+1}$,
$$c_{n+1, 0}^{\mathcal{N}}(x)= \langle f_1, \ldots, f_q, g_1, \ldots, g_q \rangle$$
$$\textrm{where } f_i= \left\{ \begin{array}{ll}
0 & \textrm{if for all } y \in N_n, \textrm{ if } c_{n, 0}^{\mathcal{N}}(y)= \alpha_i \textrm{ then } \mathcal{N} \models y \notin_n x\\
1 & \textrm{if there exists } y \in N_n, \textrm{ s.t. } c_{n, 0}^{\mathcal{N}}(y)= \alpha_i \textrm{ and } \mathcal{N} \models y \in_n x
\end{array} \right.$$
$$\textrm{and } g_i= \left\{ \begin{array}{ll}
0 & \textrm{if for all } y \in N_n, \textrm{ if } c_{n, 0}^{\mathcal{N}}(y)= \alpha_i \textrm{ then } \mathcal{N} \models y \in_n x\\
1 & \textrm{if there exists } y \in N_n \textrm{ s.t. } c_{n, 0}^{\mathcal{N}}(y)= \alpha_i \textrm{ and } \mathcal{N} \models y \notin_n x
\end{array} \right.$$

We first show that there are $\mathcal{L}_{\mathrm{TST}}$-formulae $\Phi_{n+1, 0, \beta}$, for all $\beta \in \mathcal{C}_{n+1, 0}$, that satisfy condition (II) above for the colouring $c_{n+1, 0}^{\mathcal{N}}$.
   
\begin{Lemma1} \label{Th:LiftedColouringDefinable}
For all $\beta \in \mathcal{C}_{n+1, 0}$, there is an $\mathcal{L}_{\mathrm{TST}}$-formula $\Phi_{n+1, 0, \beta}(x^{n+1})$ such that for all $x \in N_{n+1}$,
$$\mathcal{N} \models \Phi_{n+1, 0, \beta}[x] \textrm{ if and only if } c_{n+1, 0}^{\mathcal{N}}(x)= \beta.$$
\end{Lemma1}

\begin{proof}
For all $1 \leq i \leq q$, let $\Phi_{n, 0, \alpha_i}(x^{n})$ be such that for all $x \in N_n$,
$$\mathcal{N} \models \Phi_{n, 0, \alpha_i}[x] \textrm{ if and only if } c_{n, 0}^{\mathcal{N}}(x)= \alpha_i.$$
Let $\beta= \langle f_1, \ldots, f_q, g_1, \ldots, g_q \rangle \in \mathcal{C}_{n+1, 0}$. For all $1 \leq i \leq q$ and $j \in \{ 0, 1 \}$ define the $\mathcal{L}_{\mathrm{TST}}$-formula $\Theta_{i, j}^\beta(x^{n+1})$ by:
$$\Theta_{i, 0}^\beta(x^{n+1}) \textrm{ is } \left\{ \begin{array}{ll}
\forall y^n(\Phi_{n, 0, \alpha_i}(y^n) \Rightarrow y^n \notin x^{n+1}) & \textrm{if } f_i= 0\\
\exists y^n( y^n \in x^{n+1} \land \Phi_{n, 0, \alpha_i}(y^n)) & \textrm{if } f_i= 1
\end{array}\right.$$
$$\Theta_{i, 1}^\beta(x^{n+1}) \textrm{ is } \left\{ \begin{array}{ll}
\forall y^n(\Phi_{n, 0, \alpha_i}(y^n) \Rightarrow y^n \in x^{n+1}) & \textrm{if } g_i= 0\\
\exists y^n( y^n \notin x^{n+1} \land \Phi_{n, 0, \alpha_i}(y^n)) & \textrm{if } g_i= 1
\end{array}\right.$$
Define $\Phi_{n+1, 0, \beta}(x^{n+1})$ to be the $\mathcal{L}_{\mathrm{TST}}$-formula
$$\bigwedge_{1 \leq i \leq q} \bigwedge_{j \in \{0, 1\}} \Theta_{i, j}^\beta(x^{n+1}).$$
It follows from the definition of $c_{n+1, 0}^{\mathcal{N}}$ that for all $x \in N_{n+1}$,
$$\mathcal{N} \models \Phi_{n+1, 0, \beta}[x] \textrm{ if and only if } c_{n+1, 0}^{\mathcal{N}}(x)= \beta.$$
\Square
\end{proof}

We now turn to showing that $c_{n+1, 0}^{\mathcal{M}}$ and $c_{n+1, 0}^{\mathcal{N}}$ are $(2^K)^{k^\prime+2}$-similar. In order to prove this we introduce the following sets:
$$\mathrm{FOR}_n= \{ i \in [q] \mid \alpha_i \textrm{ is forbidden w.r.t. } c_{n, 0}^{\mathcal{M}} \textrm{ and } c_{n, 0}^{\mathcal{N}} \},$$
$$m\textrm{-}\mathrm{SPC}_n= \{ i \in [q] \mid \alpha_i \textrm{ is } m\textrm{-special w.r.t. } c_{n, 0}^{\mathcal{M}} \textrm{ and } c_{n, 0}^{\mathcal{N}} \} \textrm{ for } 1 \leq m < (2^K)^{k^\prime+2},$$
$$\mathrm{ABN}_n= \{ i \in [q] \mid \alpha_i \textrm{ is } (2^K)^{k^\prime+2}\textrm{-abundant w.r.t. } c_{n, 0}^{\mathcal{M}} \textrm{ and } c_{n, 0}^{\mathcal{N}} \}.$$
We classify the colours in $\mathcal{C}_{n+1, 0}$ which are forbidden, $1$-special and abundant with respect to $c_{n+1, 0}^{\mathcal{M}}$ and $c_{n+1, 0}^{\mathcal{N}}$.

\begin{Lemma1} \label{Th:ClassifyForbiddenBase}
Let $\beta \in \mathcal{C}_{n+1, 0}$ with $\beta= \langle f_1, \ldots, f_q, g_1, \ldots, g_q \rangle$. The colour $\beta$ is forbidden with respect to $c_{n+1, 0}^{\mathcal{M}}$ and $c_{n+1, 0}^{\mathcal{N}}$ if and only if either
\begin{itemize}
\item[(i)] there exists an $i \in [q]$ with $i \notin \mathrm{FOR}_n$ such that $f_i= g_i= 0$ OR,
\item[(ii)] there exists an $i \in 1\textrm{-}\mathrm{SPC}_n$ such that $f_i= g_i= 1$ OR,
\item[(iii)] there exists an $i \in \mathrm{FOR}_n$ such that $f_i= 1$ or $g_i=1$. 
\end{itemize} 
\end{Lemma1}

\begin{proof}
It is clear that if any of the conditions (i)-(iii) hold then the colour $\beta$ is forbidden. Conversely, suppose that none of the conditions (i)-(iii) hold. We need to show that $\beta$ is not forbidden with respect to $c_{n+1, 0}^{\mathcal{M}}$ and $c_{n+1, 0}^{\mathcal{N}}$. We first construct a point in $\mathcal{N}$ that is given colour $\beta$ by $c_{n+1, 0}^{\mathcal{N}}$. For all $1 \leq i \leq q$, let $\Phi_{n, 0, \alpha_i}(x^{n})$ be such that for all $x \in N_n$,
$$\mathcal{N} \models \Phi_{n, 0, \alpha_i}[x] \textrm{ if and only if } c_{n, 0}^{\mathcal{N}}(x)= \alpha_i.$$
Let $\Theta_1(x^n)$ be the $\mathcal{L}_{\mathrm{TST}}$-formula
$$\bigvee_{g_i=0} \Phi_{n, 0, \alpha_i}(x^n).$$
We work inside $\mathcal{N}$. Let $X_1= \{ x^n \mid \Theta_1(x^n) \}$. Note that comprehension ensures that $X_1$ exists. Let
$$B= \mathrm{ABN}_n \cup \bigcup_{2 \leq m < (2^K)^{k^\prime+2}} m \textrm{-}\mathrm{SPC}_n$$
and let $A= \{ i \in B \mid f_i=g_i=1 \}$. Let $\Theta_2(x^n)$ be the $\mathcal{L}_{\mathrm{TST}}$-formula
$$\bigvee_{i \in A} \Phi_{n, 0, \alpha_i}(x^n).$$
Let $X_2= \{ x^n \mid \Theta_2(x^n) \}$. Again, comprehension ensures that $X_2$ exists. For all $i \in A$, let $x_i \in N_n$ be such that $c_{n, 0}^{\mathcal{N}}(x_i)= \alpha_i$. Now, let $X= X_1 \cup (X_2 \backslash \{ x_i \mid i \in A \})$. Comprehension guarantees that $X$ exists in $\mathcal{N}$ and our construction ensures that $c_{n+1, 0}^{\mathcal{N}}(X)= \beta$. An identical construction shows that if none of the conditions (i)-(iii) hold then there is a point $X$ in $\mathcal{M}$ such that $c_{n+1, 0}^{\mathcal{M}}(X)= \beta$.     
\Square
\end{proof}

\begin{Lemma1} \label{Th:ClassifyOneSpecialBase}
Let $\beta \in \mathcal{C}_{n+1, 0}$ with $\beta= \langle f_1, \ldots, f_q, g_1, \ldots, g_q \rangle$. The colour $\beta$ is $1$-special with respect to $c_{n+1, 0}^{\mathcal{M}}$ and $c_{n+1, 0}^{\mathcal{N}}$ if and only if $\beta$ is not forbidden with respect to $c_{n+1, 0}^{\mathcal{M}}$ and $c_{n+1, 0}^{\mathcal{N}}$ and for all $i \in [q]$ with $i \notin \mathrm{FOR}_n$, $f_i= 0$ or $g_i= 0$.
\end{Lemma1}

\begin{proof}
Suppose $\beta$ is not forbidden with respect to $c_{n+1, 0}^{\mathcal{M}}$ and $c_{n+1, 0}^{\mathcal{N}}$ and for all $i \in [q]$ with $i \notin \mathrm{FOR}_n$, $f_i= 0$ or $g_i= 0$. If $x$ is a point that is given colour $\beta$ by $c_{n+1, 0}^{\mathcal{M}}$ or $c_{n+1, 0}^{\mathcal{N}}$ then $x$ is completely determined in $\mathcal{M}$ or $\mathcal{N}$ respectively. Therefore $\beta$ is $1$-special.\\ 
Conversely, suppose that $\beta$ is not forbidden and there exists an $i \in [q]$ with $i \notin \mathrm{FOR}_n$ such that $f_i= g_i= 1$. We will show that $\beta$ is not $1$-special with respect to $c_{n+1, 0}^{\mathcal{M}}$ or $c_{n+1, 0}^{\mathcal{N}}$. We first construct two distinct points of $\mathcal{N}$ that are given colour $\beta$ by $c_{n+1, 0}^{\mathcal{N}}$. For all $1 \leq i \leq q$, let $\Phi_{n, 0, \alpha_i}(x^{n})$ be such that for all $x \in N_n$,
$$\mathcal{N} \models \Phi_{n, 0, \alpha_i}[x] \textrm{ if and only if } c_{n, 0}^{\mathcal{N}}(x)= \alpha_i.$$
We work inside $\mathcal{N}$. Let $A= \{i \in [q] \mid f_i= g_i= 1 \}$. Since $\beta$ is not forbidden, for all $i \in A$, we can find $x_i, y_i \in N_n$ such that $c_{n, 0}^{\mathcal{N}}(x_i)= c_{n, 0}^{\mathcal{N}}(y_i)= \alpha_i$ and $x_i \neq y_i$. Let $\Theta_1(x^n)$ be the $\mathcal{L}_{\mathrm{TST}}$-formula
$$\bigvee_{g_i=0} \Phi_{n, 0, \alpha_i}(x^n).$$
Let $\Theta_2(x^n)$ be the $\mathcal{L}_{\mathrm{TST}}$-formula
$$\bigvee_{i \in A} \Phi_{n, 0, \alpha_i}(x^n).$$
Let $X_1= \{ x^n \mid \Theta_1(x^n) \}$ and let $X_2= \{ x^n \mid \Theta_2(x^n) \}$. Comprehension guarantees that both $X_1$ and $X_2$ exist. Let $X= X_1 \cup (X_2 \backslash \{ x_i \mid i \in A \})$ and let $Y= X_1 \cup (X_2 \backslash \{ y_i \mid i \in A \})$. Now, this construction ensures that $c_{n+1, 0}^{\mathcal{N}}(X)= c_{n+1, 0}^{\mathcal{N}}(Y)= \beta$ and $X \neq Y$. Therefore $\beta$ is not $1$-special with respect to $c_{n+1, 0}^{\mathcal{N}}$. An identical construction shows that $\beta$ is not $1$-special with respect to $c_{n+1, 0}^{\mathcal{M}}$.      
\Square
\end{proof}

\begin{Lemma1} \label{Th:ClassifyAbundantBase}
Let $\beta \in \mathcal{C}_{n+1, 0}$ with $\beta= \langle f_1, \ldots, f_q, g_1, \ldots, g_q \rangle$. If $\beta$ is not forbidden with respect to $c_{n+1, 0}^{\mathcal{M}}$ and $c_{n+1, 0}^{\mathcal{N}}$ and there exists an $i \in \mathrm{ABN}_n$ such that $f_i= g_i= 1$ then $\beta$ is $(2^K)^{k^\prime+2}$-abundant with respect to $c_{n+1, 0}^{\mathcal{M}}$ and $c_{n+1, 0}^{\mathcal{N}}$. 
\end{Lemma1}

\begin{proof}
Suppose that $\beta$ is not forbidden with respect to $c_{n+1, 0}^{\mathcal{M}}$ and $c_{n+1, 0}^{\mathcal{N}}$ and there exists an $i \in \mathrm{ABN}_n$ such that $f_i= g_i= 1$. We first construct $(2^K)^{k^\prime+2}$ distinct points in $\mathcal{N}$ that are given colour $\beta$ by $c_{n+1, 0}^{\mathcal{N}}$. For all $1 \leq i \leq q$, let $\Phi_{n, 0, \alpha_i}(x^{n})$ be such that for all $x \in N_n$,
$$\mathcal{N} \models \Phi_{n, 0, \alpha_i}[x] \textrm{ if and only if } c_{n, 0}^{\mathcal{N}}(x)= \alpha_i.$$
We work inside $\mathcal{N}$. Let $u \in \mathrm{ABN}_n$ be such that $f_u= g_u= 1$. Let $A= \{ i \in [q] \mid f_i= g_i= 1 \}$. For all $i \in A$ with $i \neq u$, let $x_i \in N_n$ be such that $c_{n, 0}^{\mathcal{N}}(x_i)= \alpha_i$. Let $y_1, \ldots, y_{(2^K)^{k^\prime+2}} \in N_n$ be such that for all $1 \leq v \leq (2^K)^{k^\prime+2}$, $c_{n, 0}^{\mathcal{N}}(y_v)= \alpha_u$ and for all $1 \leq v_1 < v_2 \leq (2^K)^{k^\prime+2}$, $y_{v_1} \neq y_{v_2}$. Let $\Theta_1(x^n)$ be the $\mathcal{L}_{\mathrm{TST}}$-formula
$$\bigvee_{g_i=0} \Phi_{n, 0, \alpha_i}(x^n).$$
Let $\Theta_2(x^n)$ be the $\mathcal{L}_{\mathrm{TST}}$-formula
$$\bigvee_{i \in A} \Phi_{n, 0, \alpha_i}(x^n).$$
Let $X_1= \{ x^n \mid \Theta_1(x^n) \}$ and let $X_2= \{ x^n \mid \Theta_2(x^n) \}$. Comprehension guarantees that $X_1$ and $X_2$ exist. For all $1 \leq v \leq (2^K)^{k^\prime+2}$, let
$$Y_v= X_1 \cup (X_2 \backslash (\{ x_i \mid i \in A \land i \neq u \} \cup \{y_v\})).$$
This construction ensures that for all $1 \leq v_1 < v_2 \leq (2^K)^{k^\prime+2}$, $Y_{v_1} \neq Y_{v_2}$ and for all $1 \leq v \leq (2^K)^{k^\prime+2}$, $c_{n+1, 0}^{\mathcal{N}}(Y_v)= \beta$. Therefore $\beta$ is $(2^K)^{k^\prime+2}$-abundant with respect to $c_{n+1, 0}^{\mathcal{N}}$. An identical construction shows that $\beta$ is $(2^K)^{k^\prime+2}$-abundant with respect to $c_{n+1, 0}^{\mathcal{M}}$.  
\Square
\end{proof}

\noindent This allows us to show that the colourings $c_{n+1, 0}^{\mathcal{M}}$ and $c_{n+1, 0}^{\mathcal{N}}$ are $(2^K)^{k^\prime+2}$-similar. 

\begin{Lemma1} \label{Th:BaseColouringsSimilar}
The colourings $c_{n+1, 0}^{\mathcal{M}}$ and $c_{n+1, 0}^{\mathcal{N}}$ are $(2^K)^{k^\prime+2}$-similar.
\end{Lemma1}

\begin{proof}
Lemma \ref{Th:ClassifyForbiddenBase} shows that for all $\beta \in \mathcal{C}_{n+1,0}$,
$$\beta \textrm{ is forbidden w.r.t. } c_{n+1, 0}^{\mathcal{M}} \textrm{ if and only if } \beta \textrm{ is forbidden w.r.t. } c_{n+1, 0}^{\mathcal{N}}.$$
Lemma \ref{Th:ClassifyOneSpecialBase} shows that for all $\beta \in \mathcal{C}_{n+1,0}$,
$$\beta \textrm{ is } 1\textrm{-special w.r.t. } c_{n+1, 0}^{\mathcal{M}} \textrm{ if and only if } \beta \textrm{ is } 1\textrm{-special w.r.t. } c_{n+1, 0}^{\mathcal{N}}.$$
Let $\beta \in \mathcal{C}_{n+1, 0}$ with $\beta= \langle f_1, \ldots, f_q, g_1, \ldots, g_q \rangle$. Lemma \ref{Th:ClassifyAbundantBase} shows that if $\beta$ is not forbidden with respect to $c_{n+1, 0}^{\mathcal{M}}$ and $c_{n+1, 0}^{\mathcal{N}}$ and there is an $i \in \mathrm{ABN}_n$ such that $f_i= g_i= 1$ then $\beta$ is $(2^K)^{k^\prime+2}$-abundant with respect to both $c_{n+1, 0}^{\mathcal{M}}$ and $c_{n+1, 0}^{\mathcal{N}}$. The remaining case is if $\beta$ is not forbidden or $1$-special and for all $i \in \mathrm{ABN}_n$, $f_i= 0$ or $g_i= 0$. Let 
$$B= \bigcup_{2 \leq m < (2^K)^{k^\prime+2}} m \textrm{-}\mathrm{SPC}.$$
In this case the number of $x \in M_{n+1}$ ($\in N_{n+1}$) with colour $\beta$ is completely determined by the number of $y \in M_{n}$ ($\in N_n$ respectively) with colour $\alpha_i$ such that $i \in B$ and $f_i= g_i= 1$. Therefore, the colourings $c_{n+1, 0}^{\mathcal{M}}$ and $c_{n+1, 0}^{\mathcal{N}}$ are $(2^K)^{k^\prime+2}$-similar.    
\Square
\end{proof}

\noindent Therefore, by induction, for all $i \in \mathbb{N}$, the colourings $c_{i, 0}^\mathcal{M}: M_i \longrightarrow \mathcal{C}_{i, 0}$ and $c_{i, 0}^\mathcal{N}: N_i \longrightarrow \mathcal{C}_{i, 0}$ are $(2^K)^{k^\prime+2}$-similar.\\
\\
\indent We now turn to defining the colour classes $\mathcal{C}_{i, j}$, and the colourings $c_{i, j}^\mathcal{M}: M_i \longrightarrow \mathcal{C}_{i, j}$ and $c_{i, j}^\mathcal{N}: N_i \longrightarrow \mathcal{C}_{i, j}$ for $1 \leq j \leq k^\prime$ and $i \in \mathbb{N}$. Let $0 \leq n < k^\prime$. Suppose that the colour classes $\mathcal{C}_{i, n}$ have been defined for all $i \in \mathbb{N}$ and that each of these colour classes has a canonical ordering. Let $j^\prime= \sum_{1 \leq m \leq n} K_m$ and suppose that $b_1^{r_1}, \ldots, b_{j^\prime}^{r_{j^\prime}} \in \mathcal{N}$ have been chosen. Moreover, suppose that for all $i \in \mathbb{N}$ and for all $\alpha \in \mathcal{C}_{i, n}$, the colourings $c_{i, n}^\mathcal{M}: M_i \longrightarrow \mathcal{C}_{i, n}$ and $c_{i, n}^\mathcal{N}: N_i \longrightarrow \mathcal{C}_{i, n}$, and the $\mathcal{L}_{\mathrm{TST}}$-formulae $\Phi_{i, n, \alpha}(x^i, \vec{z})$ have been defined with the following properties
\begin{itemize}
\item[(I')] $c_{i, n}^\mathcal{M}$ and $c_{i, n}^\mathcal{N}$ are $(2^K)^{k^\prime-n+2}$-similar, 
\item[(II')] for all $x \in N_i$,
$$\mathcal{N} \models \Phi_{i, n, \alpha}[x, b_1^{r_1}, \ldots, b_{j^\prime}^{r_{j^\prime}}] \textrm{ if and only if } c_{i, n}^\mathcal{N}(x)= \alpha.$$ 
\end{itemize}
Observe that $r_{j^\prime+1}= \ldots = r_{j^\prime+K_{n+1}}$ and let $r= r_{j^\prime+1}$. We will define the colour classes $\mathcal{C}_{i, n+1}$ and colourings $c_{i, n+1}^\mathcal{M}: M_i \longrightarrow \mathcal{C}_{i, n+1}$ and $c_{i, n+1}^\mathcal{N}: N_i \longrightarrow \mathcal{C}_{i, n+1}$ such that for all $i \in \mathbb{N}$, $c_{i, n+1}^\mathcal{M}$ and $c_{i, n+1}^\mathcal{N}$ are $(2^K)^{k^\prime-n+1}$-similar and the colouring $c_{i, n+1}^\mathcal{N}$ is definable in $\mathcal{N}$. In the process of achieving this goal we will identify points $b_{j^\prime+1}^{r}, \ldots, b_{j^\prime+K_{n+1}}^{r} \in N_{r}$.\\
\\
For all $0 \leq i < r-1$, define
$$\mathcal{C}_{i, n+1}= \mathcal{C}_{i, n},$$
$$c_{i, n+1}^\mathcal{M}= c_{i, n}^\mathcal{M},$$
$$c_{i, n+1}^\mathcal{N}= c_{i, n}^\mathcal{N}.$$
We now define the colour class $\mathcal{C}_{r-1, n+1}$, and the colourings $c_{r-1, n+1}^\mathcal{M}: M_{r-1} \longrightarrow \mathcal{C}_{r-1, n+1}$ and $c_{r-1, n+1}^\mathcal{N}: N_{r-1} \longrightarrow \mathcal{C}_{r-1, n+1}$. Let $\mathcal{C}_{r-2, n+1}= \mathcal{C}_{r-2, n}= \{ \alpha_1, \ldots, \alpha_q \}$ be obtained from the canonical ordering. Consider $a_{j^\prime+1}^{r}, \ldots, a_{j^\prime+K_{n+1}}^{r} \in M_{r}$ and use $\bar{a}_1, \ldots, \bar{a}_{K_{n+1}}$ to denote this sequence of elements. Define $\mathcal{C}_{r-1, n+1}= 2^{K_{n+1}} \times \mathcal{C}_{r-1, n}$ --- the set of all 0-1 sequences of length $K_{n+1}+2\cdot q$. Define $c_{r-1, n+1}^{\mathcal{M}}: M_{r-1} \longrightarrow \mathcal{C}_{r-1, n+1}$ such that for all $x \in M_{r-1}$,
$$c_{r-1, n+1}^{\mathcal{M}}(x)= \langle F_1, \ldots, F_{K_{n+1}}, f_1, \ldots, f_q, g_1, \ldots, g_q \rangle$$
$$\textrm{where } c_{r-1, n}^{\mathcal{M}}(x)= \langle f_1, \ldots, f_q, g_1, \ldots, g_q \rangle$$
$$\textrm{and } F_p= \left\{ \begin{array}{ll}
0 & \textrm{if } \mathcal{M} \models x \notin_{r-1} \bar{a}_p\\
1 & \textrm{if } \mathcal{M} \models x \in_{r-1} \bar{a}_p
\end{array}\right. \textrm{ for all } 1 \leq p \leq K_{n+1}.$$

\begin{Lemma1} \label{Th:ColourRefinement}
There exists $\bar{b}_1, \ldots, \bar{b}_{K_{n+1}} \in N_{r}$ such that $c_{r-1, n+1}^{\mathcal{M}}$ and the colouring $c_{r-1, n+1}^{\mathcal{N}}: N_{r-1} \longrightarrow \mathcal{C}_{r-1, n+1}$, defined such that for all $x \in N_{r-1}$,
$$c_{r-1, n+1}^{\mathcal{N}}(x)= \langle F_1, \ldots, F_{K_{n+1}}, f_1, \ldots, f_q, g_1, \ldots, g_q \rangle$$
\begin{equation} \label{Eq:RefinedNColouring}
\textrm{where } c_{r_{j^\prime+1}-1, n}^{\mathcal{N}}(x)= \langle f_1, \ldots, f_q, g_1, \ldots, g_q \rangle
\end{equation}
$$\textrm{and } F_p= \left\{ \begin{array}{ll}
0 & \textrm{if } \mathcal{N} \models x \notin_{r-1} \bar{b}_p\\
1 & \textrm{if } \mathcal{N} \models x \in_{r-1} \bar{b}_p
\end{array}\right. \textrm{ for all } 1 \leq p \leq K_{n+1},$$
are $(2^K)^{k^\prime-n+1}$-similar.  
\end{Lemma1}

\begin{proof}
Let $\mathcal{C}_{r-1, n}= \{ \alpha_1, \ldots, \alpha_{q^\prime} \}$ be obtained from the canonical ordering. For all $1 \leq i \leq q^\prime$ and for all $\sigma \in 2^{K_{n+1}}$ define $X_\sigma^i \subseteq M_{r-1}$ by
$$X_\sigma^i= \{ x \in M_{r-1} \mid (c_{r-1, n}^{\mathcal{M}}(x)= \alpha_i) \land (\forall v \in K_{n+1})(\sigma(v)=1 \iff x \in \bar{a}_v)\}.$$
Note that for all $1 \leq i \leq q^\prime$, the sets $\langle X_\sigma^i \mid \sigma \in 2^{K_{n+1}} \rangle$ partition the elements of $M_{r-1}$ that are given colour $\alpha_i$ by $c_{r-1, n}^{\mathcal{M}}$ into $2^{K_{n+1}}$ pieces. For each $1 \leq i \leq q^\prime$ choose a sequence $\langle Z_\sigma^i \mid \sigma \in 2^{K_{n+1}} \rangle$ such that for all $\sigma \in 2^{K_{n+1}}$,
\begin{itemize}
\item[(i)] $Z_\sigma^i \in N_{r}$,
\item[(ii)] for all $z \in N_{r-1}$ with $\mathcal{N} \models z \in_{r-1} Z_\sigma^i$, $c_{r-1, n}^{\mathcal{N}}(z)= \alpha_i$, 
\item[(iii)] if $|X_\sigma^i| < (2^K)^{k^\prime-n+1}$ then $|\{z \in \mathcal{N} \mid \mathcal{N} \models z \in_{r-1} Z_\sigma^i \}|= |X_\sigma^i|$,
\item[(iv)] if $|X_\sigma^i| \geq (2^K)^{k^\prime-n+1}$ then $|\{z \in \mathcal{N} \mid \mathcal{N} \models z \in_{r-1} Z_\sigma^i \}| \geq (2^K)^{k^\prime-n+1}$. 
\end{itemize}
To see that we can make this choice we work inside $\mathcal{N}$. For all $1 \leq i \leq q^\prime$, let $\Phi_{r-1, n, \alpha_i}(x^{r-1}, \vec{z})$ be such that for all $x \in N_{r-1}$,
$$\mathcal{N} \models \Phi_{r-1, n, \alpha_i}[x, b_1^{r_1}, \ldots, b_{j^\prime}^{r_{j^\prime}}] \textrm{ if and only if } c_{r-1, n}^\mathcal{N}(x)= \alpha_i.$$
For all $1 \leq i \leq q^\prime$, let $W_i= \{ x^{r-1} \mid \Phi_{r-1, n, \alpha_i}(x^{r-1}, b_1^{r_1}, \ldots, b_{j^\prime}^{r_{j^\prime}}) \}$. Comprehension ensures that the $W_i$s exist. For all $1 \leq i \leq q^\prime$ and for all $\sigma \in 2^{K_{n+1}}$, $Z_\sigma^i$ can be chosen to be a finite or cofinite subset of $W_i$. Moreover, the fact that $c_{r-1, n}^\mathcal{M}$ and $c_{r-1, n}^\mathcal{N}$ are $(2^K)^{k^\prime-n+2}$-similar ensures that for all $1 \leq i \leq q^\prime$ we can choose the sequence $\langle Z_\sigma^i \mid \sigma \in 2^{K_{n+1}} \rangle$ to satisfy condition (iii) above.\\
Now, for all $1 \leq p \leq K_{n+1}$, let $\bar{b}_p \in N_{r}$ be such that  
$$\mathcal{N} \models \bar{b}_p= \bigcup_{1 \leq i \leq q^\prime} \bigcup_{^{\sigma \in 2^{K_{n+1}}}_{\textrm{s.t. } \sigma(p)=1}} Z_\sigma^i.$$
This construction ensures that the colourings $c_{r-1, n+1}^{\mathcal{M}}$ and $c_{r-1, n+1}^{\mathcal{N}}$ define by (\ref{Eq:RefinedNColouring}) are $(2^K)^{k^\prime-n+1}$-similar.   
\Square
\end{proof}

Let $b_{j^\prime+1}^{r}, \ldots, b_{j^\prime+K_{n+1}}^{r} \in \mathcal{N}$ be the points $\bar{b}_1, \ldots, \bar{b}_{K_{n+1}}$ produced in the proof of Lemma \ref{Th:ColourRefinement} and let $c_{r-1, n+1}^{\mathcal{N}}$ be defined by (\ref{Eq:RefinedNColouring}). Therefore $c_{r-1, n+1}^{\mathcal{M}}$ and $c_{r-1, n+1}^{\mathcal{N}}$ are $(2^K)^{k^\prime-n+1}$-similar. We can immediately observe that the colouring $c_{r-1, n+1}^{\mathcal{N}}$ is definable in $\mathcal{N}$ by an $\mathcal{L}_{\mathrm{TST}}$-formula using parameters $b_{1}^{r_{1}}, \ldots, b_{j^\prime+K_{n+1}}^{r_{j^\prime+K_{n+1}}}$.

\begin{Lemma1} \label{Th:ColourRefinmentDefinable}
For all $\alpha \in \mathcal{C}_{r-1, n+1}$, there exists an $\mathcal{L}_{\mathrm{TST}}$-formula $\Phi_{r-1, n+1, \alpha}(x^{r-1}, \vec{z})$ such that for all $x \in N_{r-1}$,
$$\mathcal{N} \models \Phi_{r-1, n+1, \alpha}[x, b_{1}^{r_{1}}, \ldots, b_{j^\prime+K_{n+1}}^{r_{j^\prime+K_{n+1}}}] \textrm{ if and only if } c_{r-1, n+1}^{\mathcal{N}}(x)= \alpha.$$
\Square 
\end{Lemma1}

Let $t= \sum_{1 \leq m \leq n+1} K_m$. Lemma \ref{Th:ColourRefinement} and Lemma \ref{Th:ColourRefinmentDefinable} show that we can define colourings $c_{r-1, n+1}^{\mathcal{M}}$ and $c_{r-1, n+1}^{\mathcal{N}}$, and $\mathcal{L}_{\mathrm{TST}}$-formulae $\Phi_{r-1, n+1, \alpha}(x^{r-1}, \vec{z})$ for all $\alpha \in \mathcal{C}_{r-1, n+1}$ which satisfy the following properties:
\begin{itemize}
\item[(I'')] $c_{r-1, n+1}^{\mathcal{M}}$ and $c_{r-1, n+1}^{\mathcal{N}}$ are $(2^K)^{k^\prime-n+1}$-similar,
\item[(II'')] for all $x \in N_{r-1}$,
$$\mathcal{N} \models \Phi_{r-1, n+1, \alpha}[x, b_1^{r_1}, \ldots, b_t^{r_t}] \textrm{ if and only if } c_{r-1, n+1}^{\mathcal{N}}(x)= \alpha.$$
\end{itemize}
\indent We now turn to defining the colour classes $\mathcal{C}_{i, n+1}$, and the colourings $c_{i, n+1}^\mathcal{M}: M_i \longrightarrow \mathcal{C}_{i, n+1}$ and $c_{i, n+1}^\mathcal{N}: N_i \longrightarrow \mathcal{C}_{i, n+1}$ for all $i \geq r$. Let $i \geq r-1$. Suppose that the colour class $\mathcal{C}_{i, n+1}$ has been defined with a canonical ordering. Suppose, also, that the colourings $c_{i, n+1}^{\mathcal{M}}: M_i \longrightarrow \mathcal{C}_{i, n+1}$ and $c_{i, n+1}^{\mathcal{N}}: N_i \longrightarrow \mathcal{C}_{i, n+1}$, and the $\mathcal{L}_{\mathrm{TST}}$-formulae $\Phi_{i, n+1, \alpha}(x^{i}, \vec{z})$ have been defined and satisfy:   
\begin{itemize}
\item[(I''')] $c_{i, n+1}^{\mathcal{M}}$ and $c_{i, n+1}^{\mathcal{N}}$ are $(2^K)^{k^\prime-n+1}$-similar,
\item[(II''')] for all $x \in N_i$,
$$\mathcal{N} \models \Phi_{i, n+1, \alpha}[x, b_1^{r_1}, \ldots, b_t^{r_t}] \textrm{ if and only if } c_{i, n+1}^{\mathcal{N}}(x)= \alpha.$$
\end{itemize}
We `lift' the colour class $\mathcal{C}_{i, n+1}$ and the colourings $c_{i, n+1}^{\mathcal{M}}$ and $c_{i, n+1}^{\mathcal{N}}$ in the same way that we `lifted' the colour classes $\mathcal{C}_{i, 0}$ and the colourings $c_{i, 0}^{\mathcal{M}}$ and $c_{i, 0}^{\mathcal{N}}$ above. Let $\mathcal{C}_{i, n+1}= \{ \alpha_1, \ldots, \alpha_q \}$ be obtained from the canonical ordering. Define $\mathcal{C}_{i+1, n+1}= 2^{2 \cdot q}$--- the set of all 0-1 sequence of length $2 \cdot q$. Define $c_{i+1, n+1}^{\mathcal{M}}: M_{i+1} \longrightarrow \mathcal{C}_{i+1, n+1}$ such that for all $x \in M_{i+1}$,
$$c_{i+1, n+1}^{\mathcal{M}}(x)= \langle f_1, \ldots, f_q, g_1, \ldots, g_q \rangle$$
$$\textrm{where } f_p = \left\{ \begin{array}{ll}
0 & \textrm{if for all } y \in M_i, \textrm{ if } c_{i,n+1}^{\mathcal{M}}(y)= \alpha_p \textrm{ then } \mathcal{M} \models y \notin_i x\\
1 & \textrm{if there exists } y \in M_i \textrm{ such that } c_{i, n+1}^{\mathcal{M}}(y)= \alpha_p \textrm{ and } \mathcal{M} \models y \in_i x
\end{array}\right.$$
$$\textrm{and } g_p = \left\{ \begin{array}{ll}
0 & \textrm{if for all } y \in M_i, \textrm{ if } c_{i, n+1}^{\mathcal{M}}(y)= \alpha_p \textrm{ then } \mathcal{M} \models y \in_i x\\
1 & \textrm{if there exists } y \in M_i \textrm{ such that } c_{i, n+1}^{\mathcal{M}}(y)= \alpha_p \textrm{ and } \mathcal{M} \models y \notin_i x
\end{array}\right.$$      
Again, we define $c_{i+1, n+1}^{\mathcal{N}}$ identically. Define $c_{i+1, n+1}^{\mathcal{N}}: N_{i+1} \longrightarrow \mathcal{C}_{i+1, n+1}$ such that for all $x \in N_{i+1}$,
$$c_{i+1, n+1}^{\mathcal{N}}(x)= \langle f_1, \ldots, f_q, g_1, \ldots, g_q \rangle$$
$$\textrm{where } f_p= \left\{ \begin{array}{ll}
0 & \textrm{if for all } y \in N_i, \textrm{ if } c_{i, n+1}^{\mathcal{N}}(y)= \alpha_p \textrm{ then } \mathcal{N} \models y \notin_i x\\
1 & \textrm{if there exists } y \in N_i \textrm{ such that } c_{i, n+1}^{\mathcal{N}}(y)= \alpha_p \textrm{ and } \mathcal{N} \models y \in_i x
\end{array}\right.$$
$$\textrm{and } g_p= \left\{ \begin{array}{ll}
0 &  \textrm{if for all } y \in N_i, \textrm{ if } c_{i, n+1}^{\mathcal{N}}(y)= \alpha_p \textrm{ then } \mathcal{N} \models y \in_i x\\
1 & \textrm{if there exists } y \in N_i \textrm{ such that } c_{i, n+1}^{\mathcal{N}}(y)= \alpha_p \textrm{ and } \mathcal{N} \models y \notin_i x
\end{array}\right.$$

We first observe that there exists $\mathcal{L}_{\mathrm{TST}}$-formulae $\Phi_{i+1, n+1, \beta}(x^{i+1}, \vec{z})$ for each $\beta \in \mathcal{C}_{i+1, n+1}$ which witness the fact that the colouring $c_{i+1, n+1}^{\mathcal{N}}$ satisfies condition (II''').

\begin{Lemma1}
For all $\beta \in \mathcal{C}_{i+1, n+1}$, there is an $\mathcal{L}_{\mathrm{TST}}$-formula $\Phi_{i+1, n+1, \beta}(x^{i+1}, \vec{z})$ such that for all $x \in N_{i+1}$,
$$\mathcal{N} \models \Phi_{i+1, n+1, \beta}[x, b_1^{r_1}, \ldots, b_t^{r_t}] \textrm{ if and only if } c_{i+1, n+1}^{\mathcal{N}}(x)= \beta.$$
\end{Lemma1}

\begin{proof}
Identical to the proof Lemma \ref{Th:LiftedColouringDefinable} using the fact that $c_{i, n+1}^{\mathcal{N}}$ satisfies condition (II''').
\Square
\end{proof}

We now turn to showing that $c_{i+1, n+1}^{\mathcal{M}}$ and $c_{i+1, n+1}^{\mathcal{N}}$ are $(2^K)^{k^\prime-n+1}$-similar. To do this we prove analogues of Lemmata \ref{Th:ClassifyForbiddenBase}, \ref{Th:ClassifyOneSpecialBase} and \ref{Th:ClassifyAbundantBase}.
$$\mathrm{FOR}_i^{n+1}= \{ v \in [q] \mid \alpha_v \textrm{ is forbidden w.r.t. } c_{i, n+1}^{\mathcal{M}} \textrm{ and } c_{i, n+1}^{\mathcal{N}} \},$$
$$m\textrm{-}\mathrm{SPC}_i^{n+1}= \{ v \in [q] \mid \alpha_v \textrm{ is } m \textrm{-special w.r.t. } c_{i, n+1}^{\mathcal{M}} \textrm{ and } c_{i, n+1}^{\mathcal{N}} \} \textrm{ for } 1 \leq m < (2^K)^{k^\prime-n+1},$$
$$\mathrm{ABN}_i^{n+1}= \{ v \in [q] \mid \alpha_v \textrm{ is } (2^K)^{k^\prime-n+1} \textrm{-abundant w.r.t. } c_{i, n+1}^{\mathcal{M}} \textrm{ and } c_{i, n+1}^{\mathcal{N}} \}.$$

\begin{Lemma1} \label{Th:ClassifyForbiddenRefined}
Let $\beta \in \mathcal{C}_{i+1, n+1}$ with $\beta= \langle f_1, \ldots, f_q, g_1, \ldots, g_q \rangle$. The colour $\beta$ is forbidden with respect to $c_{i+1, n+1}^{\mathcal{M}}$ and $c_{i+1, n+1}^{\mathcal{N}}$ if and only if either
\begin{itemize}
\item[(i)] there exists a $v \in [q]$ with $v \notin \mathrm{FOR}_i^{n+1}$ such that $f_v= g_v= 0$ OR, 
\item[(ii)] there exists a $v \in 1\textrm{-}\mathrm{SPC}_i^{n+1}$ such that $f_v= g_v= 1$ OR, 
\item[(iii)] there exists a $v \in \mathrm{FOR}_i^{n+1}$ such $f_v= 1$ or $g_v=1$. 
\end{itemize}
\end{Lemma1}

\begin{proof}
Identical to the proof of Lemma \ref{Th:ClassifyForbiddenBase}.
\Square
\end{proof}

\begin{Lemma1} \label{Th:ClassifyOneSpecialRefined}
Let $\beta \in \mathcal{C}_{i+1, n+1}$ with $\beta= \langle f_1, \ldots, f_q, g_1, \ldots, g_q \rangle$. The colour $\beta$ is $1$-special with respect to $c_{i+1, n+1}^{\mathcal{M}}$ and $c_{i+1, n+1}^{\mathcal{N}}$ if and only if $\beta$ is not forbidden with respect to $c_{i+1, n+1}^{\mathcal{M}}$ and $c_{i+1, n+1}^{\mathcal{N}}$ and for all $v \in [q]$ with $v \notin \mathrm{FOR}_i^{n+1}$, $f_v=0$ or $g_v=0$. 
\end{Lemma1}

\begin{proof}
Identical to the proof of Lemma \ref{Th:ClassifyOneSpecialBase}.
\Square
\end{proof}

\begin{Lemma1} \label{Th:ClassifyAbundantRefined}
Let $\beta \in \mathcal{C}_{i+1, n+1}$ with $\beta= \langle f_1, \ldots, f_q, g_1, \ldots, g_q \rangle$. If $\beta$ is not forbidden with respect to $c_{i+1, n+1}^{\mathcal{M}}$ and $c_{i+1, n+1}^{\mathcal{N}}$ and there exists a $v \in \mathrm{ABN}_i^{n+1}$ with $f_v= g_v= 1$ then $\beta$ is $(2^K)^{k^\prime-n+1}$-abundant with respect to $c_{i+1, n+1}^{\mathcal{M}}$ and $c_{i+1, n+1}^{\mathcal{N}}$.
\end{Lemma1}

\begin{proof}
Identical to the proof of Lemma \ref{Th:ClassifyAbundantBase}.
\end{proof}

\noindent These results allow us to show that $c_{i+1, n+1}^{\mathcal{M}}$ and $c_{i+1, n+1}^{\mathcal{N}}$ are $(2^K)^{k^\prime-n+1}$-similar.

\begin{Lemma1} \label{Th:RefinedColouringsSimilar} 
The colourings $c_{i+1, n+1}^{\mathcal{M}}$ and $c_{i+1, n+1}^{\mathcal{N}}$ are $(2^K)^{k^\prime-n+1}$-similar.
\end{Lemma1}

\begin{proof}
Identical to the proof of Lemma \ref{Th:BaseColouringsSimilar} using Lemmata \ref{Th:ClassifyForbiddenRefined}, \ref{Th:ClassifyOneSpecialRefined} and \ref{Th:ClassifyAbundantRefined}.
\Square
\end{proof}

This recursion allows us to define the colour classes $\mathcal{C}_{n, k^\prime}$ and colourings $c_{n, k^\prime}^{\mathcal{M}}$ and $c_{n, k^\prime}^{\mathcal{N}}$ for all $n \in \mathbb{N}$, and elements $b_1^{r_1}, \ldots, b_1^{r_k} \in \mathcal{N}$. The above arguments show that for all $n \in \mathbb{N}$, $c_{n, k^\prime}^{\mathcal{M}}$ and $c_{n, k^\prime}^{\mathcal{N}}$ are $2^K$-similar. We have constructed the colourings $c_{n, k^\prime}^{\mathcal{M}}$ and $c_{n, k^\prime}^{\mathcal{N}}$ so as the colour assigned to a point $x \in \mathcal{M}$ (or $\mathcal{N}$) completely captures the set of quantifier-free formulae with parameters $a_1^{r_1}, \ldots, a_k^{r_k}$ (respectively $b_1^{r_1}, \ldots, b_k^{r_k}$) that are satisfied by $x$.  

\begin{Lemma1} \label{Th:ColouringsCaptureOneTypes}
Let $n \in \mathbb{N}$ and let $\theta(x_1^{r_1}, \ldots, x_k^{r_k}, x^n)$ be a quantifier-free $\mathcal{L}_{\mathrm{TST}}$-formula. If $x \in M_n$ and $y \in N_n$ are such that $c_{n, k^\prime}^{\mathcal{M}}(x)= c_{n, k^\prime}^{\mathcal{M}}(y)$ then
$$\mathcal{M} \models \theta[a_1^{r_1}, \ldots, a_k^{r_k}, x] \textrm{ if and only if } \mathcal{N} \models \theta[b_1^{r_1}, \ldots, b_k^{r_k}, y]$$
\end{Lemma1}

\begin{proof}
This follows immediately from the definition of the colourings $c_{n, k^\prime}^{\mathcal{M}}$ and $c_{n, k^\prime}^{\mathcal{N}}$.
\Square
\end{proof}

\noindent Our construction also ensures that if $x \in M_{n+1}$ (or $N_{n+1}$) then the colour assigned to $x$ by $c_{n+1, k^{\prime}}^{\mathcal{M}}$ (respectively $c_{n+1, k^\prime}^{\mathcal{N}}$) tells us, for all $\alpha \in \mathcal{C}_{n, k^\prime}$, whether there exists a point $y \in M_n$ (respectively $N_n$) such that $c_{n, k^\prime}^{\mathcal{M}}(y)= \alpha$ (respectively $c_{n, k^\prime}^{\mathcal{N}}(y)= \alpha$) and $y$ is in the relationship $\in_n$ or $\notin_n$ to $x$ in $\mathcal{M}$ (respectively $\mathcal{N}$).

\begin{Lemma1} \label{Th:KeepingTrackOfColousLemma}
Let $x \in M_{n+1}$ and $y \in N_{n+1}$, and let $\alpha \in \mathcal{C}_{n, k^\prime}$. If $c_{n+1, k^\prime}^{\mathcal{M}}(x)= c_{n+1, k^\prime}^{\mathcal{N}}(y)$ then
$$(\exists z \in M_n)(c_{n, k^\prime}^{\mathcal{M}}(z)= \alpha \land \mathcal{M} \models z \in_n x) \textrm{ if and only if } (\exists z \in N_n)(c_{n, k^\prime}^{\mathcal{N}}(z)= \alpha \land \mathcal{N} \models z \in_n y),$$
$$\textrm{and }(\exists z \in M_n)(c_{n, k^\prime}^{\mathcal{M}}(z)= \alpha \land \mathcal{M} \models z \notin_n x) \textrm{ if and only if } (\exists z \in N_n)(c_{n, k^\prime}^{\mathcal{N}}(z)= \alpha \land \mathcal{N} \models z \notin_n y).$$
\end{Lemma1}

\begin{proof}
This follows immediately from the definition of the colourings $c_{n+1, k^\prime}^{\mathcal{M}}$ and $c_{n+1, k^\prime}^{\mathcal{N}}$.
\Square
\end{proof}  

\noindent This allows us to show that an $\mathcal{L}_{\mathrm{TST}}$-sentence in the form (A) or (B) which is true $\mathcal{N}$ is also true in $\mathcal{M}$.

\begin{Theorems1} \label{Th:SentencesOfTypeBHoldInM}
Let $\phi= \forall x_1^{r_1} \cdots \forall x_k^{r_k} \exists y_1^{s} \cdots \exists y_l^{s} \theta$ be an $\mathcal{L}_{\mathrm{TST}}$-formula with $\theta$ is quantifier-free. If $\mathcal{N} \models \phi$ then $\mathcal{M} \models \phi$.
\end{Theorems1}

\begin{proof}
Suppose that $\mathcal{N} \models \phi$. Let $a_1^{r_1}, \ldots, a_k^{r_k} \in \mathcal{M}$. Using $a_1^{r_1}, \ldots, a_k^{r_k}$ and the construction we presented above we can define the colour classes $\mathcal{C}_{n, k^\prime}$ and colourings $c_{n, k^\prime}^\mathcal{M}$ and $c_{n, k^\prime}^\mathcal{N}$ for all $n \in \mathbb{N}$, and elements $b_1^{r_1}, \ldots, b_k^{r_k} \in \mathcal{N}$. The colourings $c_{n, k^\prime}^\mathcal{M}$ and $c_{n, k^\prime}^\mathcal{N}$ are $2^K$-similar and satisfy Lemma \ref{Th:ColouringsCaptureOneTypes}. Let $e_1, \ldots, e_l \in N_s$ be such that
$$\mathcal{N} \models \theta[b_1^{r_1}, \ldots, b_k^{r_k}, e_1, \ldots, e_l].$$
For all $1 \leq i \leq l$, let $d_i \in M_s$ such that $c_{s, k^\prime}^\mathcal{M}(d_i)= c_{s, k^\prime}^{\mathcal{N}}(e_i)$ and for all $1 \leq j < i$, $d_j \neq d_i$ if and only if $e_i \neq e_j$. The fact that $l < 2^K$ and $c_{s, k^\prime}^\mathcal{M}$ and $c_{s, k^\prime}^\mathcal{N}$ are $2^K$-similar ensures we can find $d_1, \ldots, d_l \in M_s$ satisfying these conditions. Now, since the variables $y_1^s, \ldots y_l^s$ all have the same type in $\theta$, the only atomic or negatomic subformulae of $\theta$ are in the form $y_i^s = y_j^s$, $y_i^s \in_s x_j^{r_j}$ if $r_j= s+1$, $x_i^{r_i} \in_{r_i} y_j^s$ if $s= r_i +1$ or $x_i^{r_i} \in_{r_i} x_j^{r_j}$ if $r_j= r_i +1$ or one of negations of these. Therefore, by Lemma \ref{Th:ColouringsCaptureOneTypes},   
$$\mathcal{M} \models \theta[a_1^{r_1}, \ldots, a_k^{r_k}, d_1, \ldots, d_l].$$
Since the $a_1^{r_1}, \ldots, a_k^{r_k} \in \mathcal{M}$ were arbitrary this shows that $\mathcal{M} \models \phi$.       
\Square
\end{proof}

\begin{Theorems1} \label{Th:SentencesOfTypeAHoldInM}
Let $\phi= \forall x_1^{r_1} \cdots \forall x_k^{r_k} \exists y_1^{s_1} \cdots \exists y_l^{s_l} \theta$ be an $\mathcal{L}_{\mathrm{TST}}$-sentence with $s_1 > \ldots > s_l$ and $\theta$ quantifier-free. If $\mathcal{N} \models \phi$ then $\mathcal{M} \models \phi$.
\end{Theorems1}

\begin{proof}
Suppose that $\mathcal{N} \models \phi$. Let $a_1^{r_1}, \ldots, a_k^{r_k} \in \mathcal{M}$. Using $a_1^{r_1}, \ldots, a_k^{r_k}$ and the construction we presented above we can define the colour classes $\mathcal{C}_{n, k^\prime}$ and colourings $c_{n, k^\prime}^\mathcal{M}$ and $c_{n, k^\prime}^\mathcal{N}$ for all $n \in \mathbb{N}$, and elements $b_1^{r_1}, \ldots, b_k^{r_k} \in \mathcal{N}$. The colourings $c_{n, k^\prime}^\mathcal{M}$ and $c_{n, k^\prime}^\mathcal{N}$ are $2^K$-similar and satisfy Lemma \ref{Th:ColouringsCaptureOneTypes}. Let $e_1^{s_1}, \ldots, e_l^{s_l} \in \mathcal{N}$ be such that
$$\mathcal{N} \models \theta[b_1^{r_1}, \ldots, b_k^{r_k}, e_1^{s_1}, \ldots, e_l^{s_l}].$$
We inductively choose $d_1^{s_1}, \ldots, d_l^{s_l} \in \mathcal{M}$. Let $d_1^{s_1} \in \mathcal{M}$ be such that $c_{s_1, k^\prime}^\mathcal{M}(d_1^{s_1})= c_{s_1, k^\prime}^{\mathcal{N}}(e_1^{s_1})$. Suppose that $1 \leq i < l$ and we have chosen $d_i^{s_i} \in \mathcal{M}$ with $c_{s_i, k^\prime}^\mathcal{M}(d_i^{s_i})= c_{s_i, k^\prime}^{\mathcal{N}}(e_i^{s_i})$. If $s_i \neq s_{i+1} + 1$ then let $d_{i+1}^{s_{i+1}} \in \mathcal{M}$ be such that $c_{s_{i+1}, k^\prime}^\mathcal{M}(d_{i+1}^{s_{i+1}})= c_{s_{i+1}, k^\prime}^{\mathcal{N}}(e_{i+1}^{s_{i+1}})$. If $s_i= s_{i+1} + 1$ and $\mathcal{N} \models e_{i+1}^{s_{i+1}} \in_{s_{i+1}} e_i^{s_i}$ then let $d_{i+1}^{s_{i+1}} \in \mathcal{M}$ be such that $c_{s_{i+1}, k^\prime}^\mathcal{M}(d_{i+1}^{s_{i+1}})= c_{s_{i+1}, k^\prime}^{\mathcal{N}}(e_{i+1}^{s_{i+1}})$ and $\mathcal{M} \models d_{i+1}^{s_{i+1}} \in_{s_{i+1}} d_i^{s_i}$. If $s_i= s_{i+1} + 1$ and $\mathcal{N} \models e_{i+1}^{s_{i+1}} \notin_{s_{i+1}} e_i^{s_i}$ then let $d_{i+1}^{s_{i+1}} \in \mathcal{M}$ be such that $c_{s_{i+1}, k^\prime}^\mathcal{M}(d_{i+1}^{s_{i+1}})= c_{s_{i+1}, k^\prime}^{\mathcal{N}}(e_{i+1}^{s_{i+1}})$ and $\mathcal{M} \models d_{i+1}^{s_{i+1}} \notin_{s_{i+1}} d_i^{s_i}$. Lemma \ref{Th:KeepingTrackOfColousLemma} and the fact that $1 < 2^K$, and $c_{s_{i+1}, k^\prime}^\mathcal{M}$ and $c_{s_{i+1}, k^\prime}^\mathcal{N}$ are $2^K$-similar ensure that we can find $d_{i+1}^{s_{i+1}} \in \mathcal{M}$ satisfying these conditions. Now, since the variables $y_1^{s_1}, \ldots y_l^{s_l}$ all have distinct types in $\theta$, the only atomic or negatomic subformulae of $\theta$ are in the form $y_{i+1}^{s_{i+1}} \in_{s_{i+1}} y_i^{s_i}$ if $s_i= s_{i+1}+1$, $y_i^{s_i} \in_{s_i} x_j^{r_j}$ if $r_j= s_i+1$, $x_i^{r_i} \in_{r_i} y_j^{s_j}$ if $s_j= r_i+1$, or $x_i^{r_i} \in_{r_i} x_j^{r_j}$ if $r_j= r_i+1$, or one of the negations of these. Therefore, by Lemma \ref{Th:ColouringsCaptureOneTypes},
$$\mathcal{M} \models \theta[a_1^{r_k}, \ldots, a_k^{r_k}, d_1^{s_1}, \ldots, d_l^{s_l}].$$
Since the $a_1^{r_1}, \ldots, a_k^{r_k} \in \mathcal{M}$ were arbitrary this shows that $\mathcal{M} \models \phi$.         
\Square
\end{proof}

\noindent Since $\mathcal{N}$ is an arbitrary model of $\mathrm{TSTI}$ and $\mathcal{M}$ is an arbitrary sufficiently large finitely generated model of $\mathrm{TST}$, Theorems \ref{Th:SentencesOfTypeBHoldInM} and \ref{Th:SentencesOfTypeAHoldInM} show that any $\mathcal{L}_{\mathrm{TST}}$-sentence in the form (A) or (B) has the finitely generated model property. Combining this with Theorem \ref{Th:ExistentialUniversalSentencesHaveFinitelyGeneratedModelProperty} shows that $\mathrm{TSTI}$ decides any sentence in the form (A) or (B).

\begin{Coroll1}
If $\phi= \forall x_1^{r_1} \cdots \forall x_k^{r_k} \exists y_1^{s_1} \cdots \exists y_l^{s_l} \theta$ is an $\mathcal{L}_{\mathrm{TST}}$-sentence with $s_1 > \ldots > s_l$ and $\theta$ quantifier free then $\mathrm{TST}$ decides $\phi$. 
\Square
\end{Coroll1}

\begin{Coroll1}
If $\phi= \forall x_1^{r_1} \cdots \forall x_k^{r_k} \exists y_1^{s} \cdots \exists y_l^{s} \theta$ is an $\mathcal{L}_{\mathrm{TST}}$-sentence with $\theta$ quantifier-free then $\mathrm{TST}$ decides $\phi$.
\Square
\end{Coroll1}

\noindent Combining these results with Theorem \ref{Th:DecidableFragmentsOfNF} shows that sentences in the form (A') or (B') are decided by $\mathrm{NF}$.

\begin{Coroll1}
If $\phi= \forall x_1 \cdots \forall x_k \exists y_1 \cdots \exists y_l \theta$ is an $\mathcal{L}$-formula with $\theta$ quantifier-free and $\sigma: \mathbf{Var}(\phi) \longrightarrow \mathbb{N}$ is a stratification of $\phi$ that assigns the same value to all of the variables $y_1, \ldots, y_l$ then $\mathrm{NF}$ decides $\phi$.
\Square
\end{Coroll1}

\begin{Coroll1}
If $\phi= \forall x_1 \cdots \forall x_k \exists y_1 \cdots \exists y_l \theta$ is an $\mathcal{L}$-formula with $\theta$ quantifier-free and $\sigma: \mathbf{Var}(\phi) \longrightarrow \mathbb{N}$ is a stratification of $\phi$ that assigns distinct values to all of the variable $y_1, \ldots, y_l$ then $\mathrm{NF}$ decides $\phi$.
\Square
\end{Coroll1}

\noindent It is interesting to note that the only use of the Axiom of Infinity in the above arguments was to ensure that the bottom type is externally infinite. Thus our arguments show that all models of $\mathrm{TST}$ with infinite bottom type agree on all sentences in the form (A) and all sentences in the form (B).

\bibliographystyle{plain}
\bibliography{decidablefragementsoftst24}

\begin{thebibliography}{1}

\bibitem{for87}
Thomas~E. Forster.
\newblock Term models for a weak set theory with a universal set.
\newblock {\em Journal of Symbolic Logic}, 52:374--387, 1987.

\bibitem{for95}
Thomas~E. Forster.
\newblock {\em Set Theory with a Universal Set: Exploring an Untyped Universe}.
\newblock Number~31 in Oxford Logic Guides. Oxford University Press, 1995.

\bibitem{hin75}
Roland Hinnion.
\newblock {\em Sur la th\'{e}orie des ensembles de Quine}.
\newblock PhD thesis, ULB, Brussels, 1975.
\newblock Translated by Thomas Forster. 2009. Available online from
  http://www.logic-center.be/Publications/Bibliotheque/hinnionthesis.pdf (last
  accessed 11/ix/2013).

\bibitem{mat01}
Adrian R.~D. Mathias.
\newblock The strength of mac lane set theory.
\newblock {\em Annals of Pure and Appied Logic}, 110(1-3):107--234, 2001.

\bibitem{qui37}
Willard v.~O. Quine.
\newblock New foundations for mathematical logic.
\newblock {\em American Mathematical Monthly}, 44:70--80, 1937.

\bibitem{rw08}
Bertrand A.~W. Russell and Alfred~N. Whitehead.
\newblock {\em Principia Mathematica}.
\newblock Cambridge University Press, 1908.

\bibitem{spe53}
Ernst~P. Specker.
\newblock The axiom of choice in quine's ``new foundations for mathematical
  logic".
\newblock {\em Proceedings of the National Academy of Sciences, U.S.A.},
  29:366--368, 1953.

\bibitem{spe62}
Ernst~P. Specker.
\newblock Typical ambiguity in logic.
\newblock In P.~Suppes E.~Nagel and A.~Tarski, editors, {\em Methodology and
  Philosophy of Science: Proceedings of the 1960 International Congress}, pages
  116--123. Stanford University Press, 1962.

\end{thebibliography}

\end{document}